\documentstyle[bezier]{article}
\catcode`\@=11
\ifcase\@ptsize
 \font\tenmsa=msam10
 \font\sevenmsa=msam7
 \font\fivemsa=msam5
 \font\tenmsb=msbm10
 \font\sevenmsb=msbm7
 \font\fivemsb=msbm5
 \font\teneu=eufm10
 \font\seveneu=eufm7
 \font\fiveeu=eufm5
 \font\tenib=cmmib10
 \font\sevenib=cmmib7
 \font\fiveib=cmmib5
\or
 \font\tenmsa=msam10 scaled \magstephalf
 \font\sevenmsa=msam7 scaled \magstephalf
 \font\fivemsa=msam5 scaled \magstephalf
 \font\tenmsb=msbm10 scaled \magstephalf
 \font\sevenmsb=msbm7 scaled \magstephalf
 \font\fivemsb=msbm5  scaled \magstephalf
 \font\teneu=eufm10  scaled \magstephalf
 \font\seveneu=eufm7  scaled \magstephalf
 \font\fiveeu=eufm5   scaled \magstephalf
 \font\tenib=cmmib10  scaled \magstephalf
 \font\sevenib=cmmib7  scaled \magstephalf
 \font\fiveib=cmmib5   scaled \magstephalf
\or
 \font\tenmsa=msam10 scaled \magstep1
 \font\sevenmsa=msam7 scaled \magstep1
 \font\fivemsa=msam5  scaled \magstep1
 \font\tenmsb=msbm10 scaled \magstep1
 \font\sevenmsb=msbm7 scaled \magstep1
 \font\fivemsb=msbm5  scaled \magstep1
 \font\teneu=eufm10   scaled \magstep1
 \font\seveneu=eufm7 scaled \magstep1
 \font\fiveeu=eufm5 scaled \magstep1
 \font\tenib=cmmib10     scaled \magstep1
 \font\sevenib=cmmib7   scaled \magstep1
 \font\fiveib=cmmib5   scaled \magstep1
\fi
\newfam\msafam
\textfont\msafam=\tenmsa  \scriptfont\msafam=\sevenmsa
  \scriptscriptfont\msafam=\fivemsa
\newfam\msbfam
\textfont\msbfam=\tenmsb  \scriptfont\msbfam=\sevenmsb
  \scriptscriptfont\msbfam=\fivemsb
\def\Bbb{\ifmmode\let\next\Bbb@\else
 \def\next{\errmessage{Use \string\Bbb\space only in math mode}}\fi\next}
\def\Bbb@#1{{\Bbb@@{#1}}}
\def\Bbb@@#1{\fam\msbfam#1}
\newfam\eufam
\textfont\eufam=\teneu  \scriptfont\eufam=\seveneu
  \scriptscriptfont\eufam=\fiveeu
\def\frak{\ifmmode\let\next\frak@\else
 \def\next{\errmessage{Use \string\frak\space only in math mode}}\fi\next}
\def\frak@#1{{\frak@@{#1}}}
\def\frak@@#1{\fam\eufam#1}
\newfam\ibfam
\textfont\ibfam=\tenib  \scriptfont\ibfam=\sevenib
  \scriptscriptfont\ibfam=\fiveib
\def\bold{\ifmmode\let\next\bold@\else
 \def\next{\errmessage{Use \string\bold\space only in math mode}}\fi\next}
\def\bold@#1{{\bold@@{#1}}}
\def\bold@@#1{\fam\ibfam#1}
\def\hexnumber@#1{\ifcase#1 0\or 1\or 2\or 3\or 4\or 5\or 6\or 7\or 8\or
 9\or A\or B\or C\or D\or E\or F\fi}
\def\newsymbolb#1#2#3#4{\mathchardef#1="#2\hexnumber@\msbfam#3#4}
\def\newsymbola#1#2#3#4{\mathchardef#1="#2\hexnumber@\msafam#3#4}
\font\fraksect=eufm10 scaled 1728
\font\fraknote=eufm8
\font\frakssect=eufm10 scaled 1440
\def\hybrid{\topmargin 0pt      \oddsidemargin 0pt
        \headheight 0pt \headsep 0pt
        \textwidth 160true mm       % US paper
        \textheight 231true mm         % US paper
        \marginparwidth 0.0in
        \parskip 0pt plus 1pt   \jot = 1.5ex}
\catcode`\@=11
\def\marginnote#1{}

\newcount\hour
\newcount\minute
\newtoks\amorpm
\hour=\time\divide\hour by60
\minute=\time{\multiply\hour by60 \global\advance\minute by-\hour}
\edef\standardtime{{\ifnum\hour<12 \global\amorpm={am}%
        \else\global\amorpm={pm}\advance\hour by-12 \fi
        \ifnum\hour=0 \hour=12 \fi
        \number\hour:\ifnum\minute<10 0\fi\number\minute\the\amorpm}}
\edef\militarytime{\number\hour:\ifnum\minute<10 0\fi\number\minute}

\def\draftlabel#1{{\@bsphack\if@filesw {\let\thepage\relax
   \xdef\@gtempa{\write\@auxout{\string
      \newlabel{#1}{{\@currentlabel}{\thepage}}}}}\@gtempa
   \if@nobreak \ifvmode\nobreak\fi\fi\fi\@esphack}
        \gdef\@eqnlabel{#1}}
\def\@eqnlabel{}
\def\@vacuum{}
\def\draftmarginnote#1{\marginpar{\raggedright\scriptsize\tt#1}}

\def\draft{\oddsidemargin -.5truein
        \def\@oddfoot{\sl preliminary draft \hfil
        \rm\thepage\hfil\sl\today\quad\militarytime}
        \let\@evenfoot\@oddfoot \overfullrule 3pt
        \let\label=\draftlabel
        \let\marginnote=\draftmarginnote
   \def\@eqnnum{(\theequation)\rlap{\kern\marginparsep\tt\@eqnlabel}%
\global\let\@eqnlabel\@vacuum}  }
\newcounter{app}
\newcounter{sapp}[app]
\def\theapp{\Alph{app}}
\newcommand{\app}[1]{
\refstepcounter{app}{\vspace{7mm}
\noindent\Large\bf Appendix
\theapp.
 \ #1 \par \vspace{5mm}}
\setcounter{equation}{0}
\def\theequation{\Alph{app}.\arabic{equation}}}

\makeatletter
\newdimen\normalarrayskip              % skip between lines
\newdimen\minarrayskip                 % minimal skip between lines
\normalarrayskip\baselineskip
\minarrayskip\jot
\newif\ifold             \oldtrue            
\def\arraymode{\ifold\relax\else\displaystyle\fi} % mode of array entries
\def\eqnumphantom{\phantom{(\theequation)}}     % right phantom in eqnarray
\def\@arrayskip{\ifold\baselineskip\z@\lineskip\z@
     \else
     \baselineskip\minarrayskip\lineskip2\minarrayskip\fi}
\def\@arrayclassz{\ifcase \@lastchclass \@acolampacol \or
\@ampacol \or \or \or \@addamp \or
   \@acolampacol \or \@firstampfalse \@acol \fi
\edef\@preamble{\@preamble
  \ifcase \@chnum
     \hfil$\relax\arraymode\@sharp$\hfil
     \or $\relax\arraymode\@sharp$\hfil
     \or \hfil$\relax\arraymode\@sharp$\fi}}
\def\@array[#1]#2{\setbox\@arstrutbox=\hbox{\vrule
     height\arraystretch \ht\strutbox
     depth\arraystretch \dp\strutbox
     width\z@}\@mkpream{#2}\edef\@preamble{\halign \noexpand\@halignto
\bgroup \tabskip\z@ \@arstrut \@preamble \tabskip\z@ \cr}%
\let\@startpbox\@@startpbox \let\@endpbox\@@endpbox
  \if #1t\vtop \else \if#1b\vbox \else \vcenter \fi\fi
  \bgroup \let\par\relax
  \let\@sharp##\let\protect\relax
  \@arrayskip\@preamble}
\def\eqnarray{\stepcounter{equation}%
              \let\@currentlabel=\theequation
              \global\@eqnswtrue
              \global\@eqcnt\z@
              \tabskip\@centering
              \let\\=\@eqncr
              $$%
 \halign to \displaywidth\bgroup
    \eqnumphantom\@eqnsel\hskip\@centering
    $\displaystyle \tabskip\z@ {##}$%
    &\global\@eqcnt\@ne \hskip 2\arraycolsep
         $\displaystyle\arraymode{##}$\hfil
    &\global\@eqcnt\tw@ \hskip 2\arraycolsep
         $\displaystyle\tabskip\z@{##}$\hfil
         \tabskip\@centering
    &{##}\tabskip\z@\cr}
\makeatother
\relax
\hybrid
\begin{document}
\def\bea{\begin{eqnarray}}
\def\eea{\end{eqnarray}}
\def\beq{\begin{equation}}          \def\bn{\beq}
\def\eeq{\end{equation}}            \def\ed{\eeq}
\def\nn{\nonumber}                  \def\g{\gamma}
\def\Uq{U_q(\widehat{\frak{sl}}_2)}
\def\Uqp{U_q(\widehat{\frak{sl}}'_2)}
\def\Uqd{U^{*}_q(\widehat{\frak{sl}}_2)}
\def\uq{U_q({sl}_2)}
\def\uqd{U^*_q({sl}_2)}
\def\slaff{\frak{sl}^\prime_2}
\def\aff{\widehat{\frak{sl}}_2}
\def\ot{\otimes}
\def\sk#1{\left({#1}\right)}
\def\id{\mbox{\rm id}}
\def\tr{\mbox{\rm tr}}
\def\tah{\mbox{\rm th}}
\def\sh{\mbox{\rm sh}}
\def\ch{\mbox{\rm ch}}
\def\ctg{\mbox{\rm ctg}}
\def\cth{\mbox{\rm cth}}
\def\tg{\mbox{\rm tg}}
\def\th{\mbox{\rm th}}
\def\qdet{\mbox{\rm q-det}}
\def\Re{{\rm Re}\,}
\def\Im{{\rm Im}\,}
\def\RR{\Bbb{R}}
\def\ZZ{\Bbb{Z}}
\def\CC{\Bbb{C}}
\def\r#1{\mbox{(}\ref{#1}\mbox{)}}
\def\d{\delta}
\def\D{\Delta}
\def\da{{\partial_\alpha}}
\let\da=p
\def\Ps{\Psi^{*}}
\def\R{{\cal R}}
\def\Ga#1{\Gamma\left(#1\right)}
\def\si#1{\sin\pi\left(#1\right)}
\def\ex#1{\exp\left(#1\right)}
\def\ep{\varepsilon}
\def\eps{\epsilon}
\def\ve{\ep}
\def\fract#1#2{{\mbox{\footnotesize $#1$}\over\mbox{\footnotesize $#2$}}}
\def\stackreb#1#2{\ \mathrel{\mathop{#1}\limits_{#2}}}
\def\res#1{\stackreb{\mbox{\rm res}}{#1}}
\def\lim#1{\stackreb{\mbox{\rm lim}}{#1}}
\def\Res#1{\stackreb{\mbox{\rm Res}}{#1}}
\let\dis=\displaystyle
\def\ee{{\rm e}}
\def\D{\Delta}
\def\H{{\cal H}}
\renewcommand{\theequation}{{\thesection}.{\arabic{equation}}}
\def\Y-{\widehat{Y}^-}
\def\DYsect{\widehat{DY(\hbox{\fraksect sl}_2)}}
\def\Ael{{\cal A}_{\tih,\eta}(\widehat{\frak{sl}_2})}
\def\Aelx{{\cal A}_{\tih,\xi}(\widehat{\frak{sl}_2})}
\def\Aelsect{{\cal A}_{\tih,\eta}(\widehat{\hbox{\fraksect sl}_2})}
\def\Aelnote{{\cal A}_{\tih,\eta}(\widehat{\hbox{\fraknote sl}_2})}
\def\Aelssect{{\cal A}_{\tih,\eta}(\widehat{\hbox{\frakssect sl}_2})}
\def\Apq{{\cal A}_{q,p}(\widehat{\frak{sl}_2})}
\def\Apqsect{{\cal A}_{q,p}(\widehat{\hbox{\fraksect sl}_2})}
\def\DY{\widehat{DY(\frak{sl}_2)}}
\def\Yd{\DY}
\def\Ydd{\DY}
\let\z=z
\let\b=z
\def\u{{u}}
\def\v{{v}}
\def\g{\gamma}
\def\la{\lambda}
\let\hsp=\qquad
\def\he{{\hat e}}
\def\hf{{\hat f}}
\def\hh{{\hat t}}
\def\ha{{\hat a}}
\def\hb{{\hat b}}
\def\hk{{\hat t}}
\def\hkp{{\hat t}}
\def\hhh{{\hat h}}
\def\Ev{{\cal E}v}
\def\vvv{\overline{\varphi}}
\def\vac{|\mbox{vac}\rangle}
\def\lvac{\langle\mbox{vac}|}
\def\vpint{-\!\!\!\!\!\!\!\int_{-\infty}^\infty}
\def\vpinto{-\!\!\!\!\!\!\!\int_{0}^\infty}
\def\vpintss{{\vpint\cdots\ \vpint}}
\def\vpintoo{{\vpinto\cdots\ \vpinto}}
\def\intt{\int_{-\infty}^\infty}
\def\vpints{-\!\!\!\!\!\!\!\int_{-\infty}^{(\la-\mu)/2}}
\def\la{\lambda}
\def\tih{{\hbar}}
\def\FDY{F\left[\DY\right]}
\def\FDYsect{F\left[\DYsect\right]}
\def\stackupb#1#2#3{\ \mathrel{\mathop{#1}\limits_{#2}^{#3}}}
\def\feq#1#2{\stackupb{\ravnodots}{#1}{#2}}
\def\cint{\int_\infty^{0+}}
\def\nint{\int^{+\infty}_{0}}
\def\mint{\int_{-\infty}^{0}}
\let\rvac=\vac
%%%%%%%%%%%%%%%%%%%%%%%%%%%%%%%%%%%%%%%%%%%%%%%%%%%%%%%%%%%%%%%%%%%%%%
%%%%%%%%%%%%%%%%%   End of Personal difinitions    %%%%%%%%%%%%%%%%%%%
%%%%%%%%%%%%%%%%%%%%%%%%%%%%%%%%%%%%%%%%%%%%%%%%%%%%%%%%%%%%%%%%%%%%%%
\begin{titlepage}
\begin{center}
\hfill ITEP-TH-51/96\\
%\hfill JINR-95-??\\
\hfill q-alg/9702002\\
\bigskip\bigskip
{\Large\bf Elliptic Algebra $\Apqsect$ in the Scaling Limit}\\
\bigskip
\bigskip
{\large S. Khoroshkin\footnote{E-mail: khoroshkin@vitep1.itep.ru},
D. Lebedev\footnote{E-mail: lebedev@vitep1.itep.ru}}\\
\medskip
{\it Institute of Theoretical \& Experimental Physics\\
117259 Moscow, Russia}\\
\bigskip
{\large S. Pakuliak}\footnote{E-mail: pakuliak@thsun1.jinr.dubna.su}\\
\bigskip
{\it Bogoliubov Laboratory of Theoretical Physics, JINR\\
141980 Dubna, Moscow region, Russia}\\
\bigskip
\bigskip
\bigskip
%{Revised \today}
\end{center}
\begin{abstract}
The scaling limit $\Ael$ of the elliptic  algebra $\Apq$ is investigated.
The limiting algebra is defined in terms of a continuous family of 
 generators being Fourier harmonics of Gauss coordinates of the $L$-operator.
 Ding-Frenkel isomorphism between $L$-operator's and current descriptions of
the algebra $\Ael$ is established and is identified with the Riemann problem 
 on a strip. The representations, coalgebraic structure and intertwining 
 operators of the algebra are studied.
\end{abstract}
\end{titlepage}
\clearpage
\newpage

\setcounter{section}{-1}
\setcounter{equation}{0}
\setcounter{footnote}{0}

\section{Introduction}
This  paper is devoted to the  investigation of the infinite-dimensional
 algebra $\Ael$ which is supposed to be the algebra of symmetries in 
 integrable models of quantum field theories. Our work arose 
from the attempts to understand the mathematical background of the results
 by
S.~Lukyanov \cite{L} and to combine his methods with the
group-theoretical approach to quantum integrable models
 developed in \cite{JMbook,FIJKMY}
 and with the Yangian calculations in 
 \cite{KLP1}.

The investigation of the symmetries in 
quantum integrable models of the two-dimensional field theory
started in  \cite{BKKW,Lu''} 
resulted in the form-factor (bootstrap) approach to these
models developed in the most completed form in the works
by F.A.~Smirnov \cite{S1}. This approach was not addressed
to investigation of dynamical symmetries in the model but 
 to computation of
 certain final objects of the theory, form-factors
of the local operators and correlation
functions  of the local operators. 
It was observed in the papers \cite{B,Sm,BL} that 
the mathematical structures underlying  the success of the
bootstrap approach   in the massive integrable models 
are related to the representation theory of infinite-dimensional Hopf
algebras.
The
dynamical symmetries in massive two-dimensional field theories
was investigated in \cite{L}  
in the framework of Zamolodchikov-Faddeev operators
\cite{ZZ,F}.

The elliptic algebra  $\Apq$ was proposed  in the works 
\cite{FIJKMY} as an algebra of symmetries
 for the eight-vertex lattice integrable model.
The algebra $\Apq$ was formulated in the framework of the ``$RLL$'' approach
\cite{FRT} in terms of the symbols $L^\pm_{\ep\ep',n}$ ($n\in\ZZ$,
$\ep,\ep'=\pm$, $\ep\ep'=(-1)^n$) gathered using the spectral parameter
$\zeta$ into
$2\times2$ matrices $L^\pm(\zeta)$ and   the central element $c$.
The generating series $L^\pm(\zeta)$ satisfy the defining relations:
\bea
R^\pm_{12}(\zeta_1/\zeta_2)L^\pm_1(\zeta_1)L^\pm_2(\zeta_2)
&=& L^\pm_2(\zeta_2)L^\pm_1(\zeta_1)R^{*\pm}_{12}(\zeta_1/\zeta_2),\nn\\
R^+_{12}(q^{c/2}\zeta_1/\zeta_2)L^+_1(\zeta_1)L^-_2(\zeta_2)
&=& L^-_2(\zeta_2)L^+_1(\zeta_1)R^{*+}_{12}(q^{-c/2}\zeta_1/\zeta_2),\nn\\
q^{c/2}&=&L^+_{++}(q^{-1}\zeta)L^+_{--}(\zeta) -
   L^+_{-+}(q^{-1}\zeta)    L^+_{+-}(\zeta),\nn\\
L^-_{\ep\ep'}(q^{-1}\zeta)&=& \ep\ep' L^+_{-\ep,-\ep'}(p^{1/2}q^{-c/2}\zeta),
\label{0.1}
\eea
where
\beq
R^\pm(\zeta)=q^{\mp1/2}\zeta\left[
{(q^2\zeta^{-2};q^4)_\infty (q^2\zeta^{2};q^4)_\infty
\over
(q^4\zeta^{\mp2};q^4)_\infty (\zeta^{\pm2};q^4)_\infty }\right]^{\pm1}
R(\zeta),
\label{Baxter}
\eeq
and $R(\zeta)=R(\zeta;p^{1/2},q^{1/2})$ is the Baxter elliptic $R$ matrix
normalized to satisfy the unitarity and crossing symmetry relations
\cite{Ba}, $R^{*\pm}=R(\zeta;{p^*}^{1/2},q^{1/2})$ and $p^*=pq^{-2c}$. 
Unfortunately, there is no description of
 infinite-dimensional representations of $\Apq$ in terms of free fields.

In this paper we are going to investigate the scaling
limit of this algebra when $q,p\to1$. We call this algebra
$\Ael$.
 Let us note  that although the algebra $\Ael$
is constructed by means of the trigonometric $R$-matrix, it
 is quite different from the quantum affine algebra which
is a degeneration of the elliptic algebra
$\Apq$ when $p=0$. The algebra $\Ael$, written in integral relations
 for usual commutators and anticommutators,
 conserve many principal  properties of the elliptic algebra. For example,
 it possesses an evaluation homomorphism onto a degenerated 
 Sklyanin algebra \cite{Sk} for zero central  charge and has no 
 imbeddings of finite\-dimensional quantum groups for $c\neq 0$.
 But due to the more simple structure of $R$-matrix comparing with 
 the elliptic 
 case  a more detailed study of its algebraical structure is possible.

 In particular, one of our achievements is  the currents  description
of the algebra $\Ael$, which is equivalent to the factorization of 
the quantum determinant 
in the $L$-operator approach. This allows us to make a more detailed 
investigation of
the representation theory of the algebra $\Ael$. Starting from basic
 representation of $\Ael$ in a Fock space we reconstruct precisely 
the  Zamolodchikov--Faddeev algebra, described in \cite{ZZ,L,JKM}. 
We carry out
 this reconstruction from the analysys of the Hopf structure of $\Ael$.

The distinguished feature of the algebra $\Ael$ is a presence of 
analysis in its
 description. The formal generators of the algebra are Fourier harmonics of
 the currents labeled by real numbers, and the elements of the algebra are
 integrals over generators with coefficients being functions with certain
 conditions on their analyticity and on their asymptotical behaviour. 

The paper is organized as follows.
In the first section  we give a description of the algebra $\Ael$ in terms of
 formal generators being Fourier harmonics of the currents. 
The relations for the formal generators are given in a simple integral form.
We assign a precise meaning to the elements of the algebra as to certain
integrals over generators and show that the relations are correctly defined in
corresponding vector spaces. Moreover, we show  that the quadratic integral
 relations could be interpreted as ordering rules with polylogarithmic 
 coefficients for monomials  composed from  generators of the algebra. In this
section we suppose that central charge is not equal to zero. In the next
section we develop the formalism of $L$-operators for $\Ael$. We show that 
the $L$-operators $L^{\pm}(u)$, satisfying the standard relations 
\cite {FRT} with
$R$-matrices being scaling limits of those from \cite{Ba}, 
 admit the Gauss decomposition. We write down 
 relations for the Gauss coordinates and identify them with generating
functions for the generators of $\Ael$ described in the previous section.
 Looking to the rational limit $\eta \to 0$ we find a double of the
 Yangian but in  a  presentation different from \cite{K,IK}.
 We describe also the coalgebraic structure of $\Ael$. The comultiplication
rule looks standard in terms of the $L$-operators, it is compatible with the
defining relations, but it sends now the initial algebra into a tensor product
 of two different algebras which differ by the value of the parameter $\eta$.
 Nevertheless it is sufficient for the definition of the
intertwining operators.
 We call this structure a Hopf family of algebras.

Section 3 is devoted to the description of algebra $\Ael$ for $c = 0$. We
treat this case as a limit of $\Ael$ when $c$ tends to zero. The limit is not
trivial, one should look carefully to the asymptotics of the currents in the 
 limit in order to define correct generators for $c = 0$. We describe 
 finite-dimensional representations and the evaluation homomorphism onto the 
 degenerated Sklyanin algebra, which is isomorphic in this case to $U_q(sl_2)$
 with $|q|=1$. 
In the next section we complete Ding-Frenkel isomorphism \cite{FD} and present
 a description of the algebra in terms of total currents. We show that 
 Ding-Frenkel formulas are equivalent in our case to Sokhotsky-Plemely's
formulas for the Riemann problem on a strip. The relation \r{0.1} for $L^\pm$
 operators is also natural in the framework of the Riemann problem.

The last two sections are devoted to the study of the basic representation 
of the algebra $\Ael$ in a Fock space. The representation of the corresponding
 Zamolodchikov-Faddeev algebra in this space was recently described in 
 \cite{JKM,Ko}. We start from a bosonization of
 the total currents for $\Ael$ and then identify Zamolodchikov-Faddeev
 algebra with the algebra of 
type I and type II twisted intertining operators. The twisting
 means a presence of a certain involution in the definition of the 
intertwining 
 operators. The twisting comes from the lack of zero mode operator 
 $(-1)^p$ in the continious models. There is no motivation to introduce this
operator in our case since, to the contrary to discrete models, we have 
the unique
 level one module. As a consequence, Zamolodchikov-Faddeev operators commute
 by means of an $R$ matrix \cite{ZZ,L}
which differs from the one used in the description of 
 $\Ael$  by certain signes. We check also the correspondence
of the Miki's formulas \cite{M} to the $L$-operator description of the  basic
representation of $\Ael$. Note also that the notions of a Fock space and of 
 vertex operators for continuous free boson field require special 
 analytical definition which we suggest in the last section.

\bigskip

\noindent {\sl
 We use the opportunity to thank A.~Belavin, M.~Jimbo,
S.~Lukyanov and  T.~Miwa 
for the stimulating discussions and
for their remarkable works which have been the source of our
inspiration. We would like to acknowledge also the useful discussions
with  A.~Gerasimov and Ya.~Pugai}.

\setcounter{equation}{0}
\section{Algebra $\Aelsect$ ($c\neq0$)}

\subsection{The definition}

For $\la\in\RR$ we consider the family
of symbols
$\he_\la$, $\hf_\la$, $\hh_\la$ and  $c$
of the formal algebra
which satisfy the commutation relations:
\bea
{[}c,\mbox{everything}{]}&=&0\ ,\label{c-ever}\\
{[}\he_\la,\hf_\mu{]}&=&
\sh\left(\fract{\la}{2\eta}+\fract{\mu}{2\eta'}\right)
\hh_{\la+\mu}
\ ,\label{he-hf}\\
{[}\hh_\la,\he_\mu{]}&=&
{\tg\,\pi\eta\tih\over 2 \pi\eta}
\vpint d\tau\
\sh\left(\fract{\tau}{2\eta}\right)^{-1}
\{\hh_{\la+\tau},\he_{\mu-\tau}\}\ ,
\label{hh-he}\\
{[}\hh_\la,\hf_\mu{]}&=&
-{\tg\,\pi\eta'\tih\over 2 \pi\eta'}
\vpint d\tau\
\sh\left(\fract{\tau}{2\eta'}\right)^{-1}
\{\hh_{\la+\tau},\hf_{\mu-\tau}\},
\label{hh-hf}\\
{[}\he_\la,\he_\mu{]}&=&
{\tg\,\pi\eta\tih\over 2 \pi\eta}
\vpint d\tau\ \cth\left(\fract{\tau}{2\eta}\right)
\{\he_{\la+\tau},\he_{\mu-\tau}\}\ ,
\label{he-he}\\
{[}\hf_\la,\hf_\mu{]}&=&
-{\tg\,\pi\eta'\tih\over 2 \pi\eta'}
\vpint d\tau\ \cth\left(\fract{\tau}{2\eta'}\right)
\{\hf_{\la+\tau},\hf_{\mu-\tau}\}\ ,
\label{hf-hf}\\
{[}\hh_\la,\hh_\mu{]}&=&
\vpint d\tau\ \kappa(\tau)
\{\hh_{\la+\tau},\hh_{\mu-\tau}\}\ ,
\label{hh-hh}
\eea
where 
the real odd function $\kappa(\tau)$ is given by the Fourier transform
$$%\beq
\kappa(\tau)={1\over2\pi} \intt du\  \ee^{-i\tau u}
{\cth(\pi\eta'u)\th(i\pi\eta'\tih)-
\cth(\pi\eta u)\th(i\pi\eta\tih)\over
1-\cth(\pi\eta'u)\cth(\pi\eta u)
\th(i\pi\eta'\tih)\th(i\pi\eta\tih)}\ ,
$$%\label{g-funk}\eeq
$\tih$ is a deformation parameter, $\{a,b\}$ means $ab+ba$,
$\eta>0$ and the 
parameters $\eta$ and $\eta'$ are related through the central element $c$:
$$%\beq
{1\over\eta'}-{1\over\eta}=\tih c\ ,\qquad \tih c>0\  .
$$%\label{eta-rel}\eeq
The last inequality means that in the
  representations which we consider  the
central element
$c$ is equal to some number such that 
$\tih c>0$  and  we identify $c$
with this number.
The case $c=0$ requires a special
treatment and will be considered in the next section.

%This restriction \r{eta-rel} corresponds to the restriction
% $|q|<1$ 
%for  highest weight representations of the quantum affine algebra
%$U_q(\widehat{\frak{sl}}_2)$.

Let us consider the vector space $\overline{{\cal A}}$
formed by the formal integrals of the type
\beq
\intt \prod_k d\la_k \prod_i d\mu_i \prod_j d\nu_j\
\phi(\{\la_k\};\{\mu_i\};\{\nu_j\})
P(\{\he_{\la_k}\};\{\hf_{\mu_i}\};\{\hh_{\nu_j}\})\ ,
\label{elements}
\eeq
where $\phi(\{\la_k\};\{\mu_i\};\{\nu_j\})$
 is the  $\CC$-number function of real
variables $\la_k$, $\mu_i$ and $\nu_j$ which
satisfy  the conditions
of analyticity:
\bea
&\phi(\{\la_k\};\{\mu_i\};\{\nu_j\})
&\quad \mbox{is analytical in the strip}\quad
-\pi\eta<\Im\la_k<\pi\eta\quad \forall\ \ \la_k\ ,\nn\\
&\phi(\{\la_k\};\{\mu_i\};\{\nu_j\})
&\quad \mbox{is analytical in the strip}\quad
-\pi\eta'<\Im\mu_i<\pi\eta'\quad \forall\ \ \mu_i\ ,\nn\\
&\phi(\{\la_k\};\{\mu_i\};\{\nu_j\})
&\quad \mbox{is analytical  in the strip}\quad
-\pi\eta'<\Im\nu_j<\pi\eta'\quad \forall\ \ \nu_j\ ,\nn
\eea
and conditions on the asymptotics when $\Re\la_k$, $\Re\mu_i$,
$\Re\nu_j\to\pm\infty$:
\bea
\phi(\{\la_k\};\{\mu_i\};\{\nu_j\})&<& C\ee^{-\alpha|\Re\la_k|}\ ,\nn\\
\phi(\{\la_k\};\{\mu_i\};\{\nu_j\})
&<& C\ee^{-(\beta+\tih c/2)|\Re\mu_i|}\ ,\nn\\
\phi(\{\la_k\};\{\mu_i\};\{\nu_j\})
&<& C\ee^{-(\gamma-1/2\eta)|\Re\nu_j|}\ ,\nn
\eea
for some real positive $\alpha$, $\beta$, $\gamma$. The notation
 $$
P(\{\he_{\la_k}\};\{\hf_{\mu_i}\};\{\hh_{\nu_j}\})
$$
means
monomial which is a product
of the formal generators $\he_\la$, $\hf_\mu$ and
$\hh_\nu$ in some order.

The space $\overline{{\cal A}}$ has a natural structure of free (topological)
 algebra. 

By definition
the algebra $\Ael$ is identified with $\overline{{\cal A}}$
factorized by the ideal generated by the commutation relations
\r{he-hf}--\r{hh-hh} which
 can be treated as equalities in
the vector space $\overline{{\cal A}}$.

The correctness of the definition of the algebra $\Ael$ follows from
the Lemma 1 and the properties of the kernels of the integral
transforms which enter in the r.h.s. of the commutation
relations \r{he-hf}--\r{hh-hh}. These relations make also possible
to write the monomials in \r{elements} in the ordered form (see the next
 subsection).
\smallskip

\noindent
{\bf Lemma 1.} {\it For two functions $a(\la)$ and $b(\la)$ which are
analytical in the strips $-\alpha_1<\Im\la<\alpha_2$,
$-\beta_1<\Im\la<\beta_2$ respectively
for $\alpha_1$,$\alpha_2$, $\beta_1$ 
$\beta_2 >0$  and have exponentially decreasing
asymptotics when $\Re\la\to\pm\infty$ the convolution}
$$%\beq
(a\star b)(\la)=
\intt \d\tau\ a(\tau)\, b(\la-\tau)
$$%\label{convolution}\eeq  
{\it
is analytical function of $\la$ in the strip
$-\alpha_1-\beta_1<\Im\la<\alpha_2+\beta_2$,
and also have exponentially decreasing asymptotics at
$\Re\la\to\pm\infty$.}
\smallskip

In the sequel
we will
need following involution  of the algebra $\Ael$:
\beq
\iota\sk{\he_\la}=-\he_\la,\quad
\iota\sk{\hf_\la}=-\hf_\la,\quad
\iota\sk{\hh_\la}=\hh_\la\ .\label{auto}
\eeq

\subsection{Commutation relations as ordering rules}

Assign the meaning to the commutation relations
\r{he-hf}--\r{hh-hh}. One should understand them as the rules to
express the product of the formal generators in the form
$$
\he_{\la_1}\he_{\la_2}\ldots
\hf_{\la_1}\hf_{\la_2}\ldots
\hh_{\la_1}\hh_{\la_2}\ldots
$$
and then as a rule to express the product, say,
$\he_\la\he_\mu$ for
$\la>\mu$ through the products        $\he_{\tau_1}\he_{\tau_2}$
for $\tau_1<\tau_2$.

It is clear how to use the commutation relation \r{he-hf} to order
the products. Let us explain the use of the rest of the commutation relations.
The commutation relations of the type \r{hh-he} and \r{hh-hf}
$$
{[}\ha_\la,\hb_\mu{]}=
\vpint d\tau\ \varphi(\tau)
\{\ha_{\la+\tau},\hb_{\mu-\tau}\}\ ,
$$
can be rewritten in the form of the ordering  rules which 
structure coefficients are composed from 
 iterated integrals of the po\-ly\-lo\-ga\-riph\-mic type:
\bea
\ha_\la\hb_\mu&=&\hb_\mu\ha_\la   +
2\vpint d\tau\
\hb_{\mu-\tau}\ha_{\la+\tau}\times\nn\\
&\times&\left(
\varphi(\tau)
+\sum_{n=1}^\infty
\vpintss d\tau_1\ldots d\tau_n\  \varphi(\tau-\tau_1)
\prod_{k=1}^{n-1} \varphi(\tau_k-\tau_{k+1}) \varphi(\tau_n)
\right),\label{ha-hb}
\eea
where
the function $\varphi(\tau)$ is
$${\th(\pi\eta\tih)\over
2\pi\eta\sh(\tau/2\eta)}\quad
\mbox{for \r{hh-he}}\quad \mbox{and}\quad
-{\th(\pi\eta'\tih)\over
2\pi\eta'\sh(\tau/2\eta')}\quad
\mbox{for \r{hh-hf}.}$$

The operators of the same type can be ordered according to the
ordering of the indeces. Fix $\la>\mu$.
Then iterating
the commutation relations of the type \r{he-he}, \r{hf-hf} or
\r{hh-hh} we obtain:
\bea
\ha_\la\ha_\mu&=&\ha_\mu\ha_\la + 2
\vpinto d\tau\
\ha_{{\la+\mu\over2}-\tau}\ha_{{\la+\mu\over2}+\tau}\times\nn\\
&\times&\!\!\!\!\!
\left(
\vvv\left(\fract{\la-\mu}{2};\tau\right) + \sum_{n=1}^\infty
\vpintoo d\tau_1\ldots d\tau_n\  \vvv\left(\fract{\la-\mu}{2};\tau_1\right)
\prod_{k=1}^{n-1} \vvv(\tau_k;\tau_{k+1}) \vvv(\tau_n;\tau)
\right),\label{ha-ha}
\eea
where
$$
\vvv(\tau;\tau')=\varphi(\tau-\tau')+\varphi(\tau+\tau')\
$$
and
the function $\varphi(\tau)$ is
$$-{\th(\pi\eta\tih)\over
2\pi\eta\th(\tau/2\eta)}\quad
\mbox{for \r{he-he},}\quad
{\th(\pi\eta'\tih)\over
2\pi\eta'\th(\tau/2\eta')}\quad
\mbox{for \r{hf-hf}}$$
and $-\kappa(\tau)$ for \r{hh-hh}.
\smallskip

\noindent
{\bf Conjecture 2.} {\it The series in}
\r{ha-hb} and \r{ha-ha} {\it are convergent}.
\smallskip

There are few remarks in favour if this conjecture. First, we see
that if we consider the deformation parameter $\tih$ small than
these series are series with respect to powers of the small
parameter. Second, in the Yangian limit when $\eta\to0$ the
series in \r{ha-ha}
can be summed up to obtain the function $\tih\ee^{\tih\tau}$.

\setcounter{equation}{0}
\section{$L$-Operator Realization of the Algebra $\Aelsect$}

\subsection{Gauss coordinates of the $L$-operator}

Fix the following $R$-matrix \cite{ZZ,L}:
\begin{eqnarray}
R^+(u,\eta)&=&\tau^+(u)  R (u,\eta),\quad
R(u,\eta)\ =\ r (u,\eta)\overline R (u,\eta)\ ,
\nn\\
\overline R (u,\eta)&=&
\left(
\begin{array}{cccc}
1&0&0&0\\  0&b (u,\eta)&c (u,\eta)&0\\
0&c (u,\eta)&b (u,\eta)&0\\  0&0&0&1
\end{array}
\right)      \ ,     \nn\\
r (z,\eta)&=&{\Ga{\tih\eta}\Ga{1+i\eta u}\over
        \Ga{\tih\eta+i\eta u}}
\prod_{p=1}^\infty
{R _p(u,\eta) R _p(i\tih- u,\eta)
\over
R _p(0,\eta) R _p(i\tih,\eta)}\ , \nn\\
R _p(u,\eta)&=&{\Ga{2p\tih\eta+i\eta u}
               \Ga{1+2p\tih\eta+i\eta u}\over
        \Ga{(2p+1)\tih\eta+i\eta u}
        \Ga{1+(2p-1)\tih\eta+i\eta u}}\ ,\nn\\
b(u,\eta)\ &=&\
{\sh\,\pi\eta u\over\sh\,\pi\eta(u-i\tih)}\ ,\quad
c(u,\eta)\ =\
{-\sh\,i\pi\eta\tih \over\sh\,\pi\eta(u-i\tih)}\ ,\quad
\tau^+(u)=\cth\sk{{\pi u\over2\tih}} ,   \label{R-mat}
\end{eqnarray}
where $u$ is a spectral parameter.

Let
\beq
L(u)=\left(\begin{array}{cc}
L_{++}(u)&L_{+-}(u)\\ L_{-+}(u)&L_{--}(u)
\end{array}\right)
\label{L-op}
\eeq
be a quantum $L$-operator which matrix elements are treated as generating
functions for the elements of the algebra given by the
commutation relations:
\begin{eqnarray}
R^+(u_1-u_2,\eta')L_1(u_1,\eta)L_2(u_2,\eta)&=&
L_2(u_2,\eta)L_1(u_1,\eta) R^+(u_1-u_2,\eta) \  , \label{RLL-univ}\\
\qdet L(u)&=&1\ .
\label{qdet=1}
\eea
The quantum determinant of the $L$-operator
is given by
\bea
\qdet L(z)
= L_{++}(z-i\tih)L_{--}(z)-L_{+-}(z-i\tih)L_{-+}(z)\ .%\nn\\
%&=& -D(z-i\tih)A(z)-C(z-i\tih)B(z)\nn\\
%&=& -D(z)A(z-i\tih)-B(z)C(z-i\tih)\nn\\
%&=& -A(z)D(z-i\tih)-C(z)B(z-i\tih)
\label{qdet}
\eea

Let
\beq
L(u)
=\left(\begin{array}{cc} 1& f(u)\\0&1\end{array}\right)
\left(\begin{array}{cc}  k_1(u)&0\\ 0&k_2(u) \end{array}\right)
\left(\begin{array}{cc} 1&0\\ e(u)&1\end{array}\right)\ ,
\label{GL-univ}
\eeq
be the Gauss decomposition of the $L$-operator \r{L-op}. One can deduce from
 (\ref{RLL-univ}), (\ref{qdet=1}) that
$$
k_1(u)=(k_2(u+i\tih))^{-1}.
$$
Let 
$$%\beq
h(u)=k_2\left(u\right)^{-1} k_1\left(u\right)\ , \qquad
\tilde h(u)=k_1\left(u\right)  k_2\left(u\right)^{-1}
\ . 
$$%\label{h-new}\eeq
Then, due to (\ref{RLL-univ}), (\ref{qdet=1})
\beq
\tilde h(u)= {\eta\over\eta'}
{\sin\,\pi\eta'\tih\over\sin\,\pi\eta\tih} h(u)\ .\label{til-h-new} 
\eeq
 We have the following 
\smallskip

\noindent
{\bf Proposition 3.} {\it The Gauss coordinates e(u), f(u) and h(u) of
 the $L$-operator} \r{L-op} {\it satisfy the following commutation relations 
} ($u=u_1-u_2$)
\bea
e(u_1)f(u_2)-f(u_2)e(u_1)
&=&{\sh\,i\pi\eta'\tih\over\sh\,\pi\eta'u}h(u_1)-
{\sh\,i\pi\eta\tih\over\sh\,\pi\eta u}\tilde h(u_2),\label{ef}\\
\sh\,\pi\eta(u+i\tih)h(u_1)e(u_2)-
\sh\,\pi\eta(u-i\tih)e(u_2)h(u_1)
&=&\sh(i\pi\eta\tih)
\{h(u_1),e(u_1)\},\label{he}\\
\sh\,\pi\eta'(u-i\tih)h(u_1)f(u_2)-
\sh\,\pi\eta'(u+i\tih)f(u_2)h(u_1)
&=&-\sh(i\pi\eta'\tih)
\{h(u_1),f(u_1)\},\label{hf}\\
\sh\,\pi\eta(u+i\tih)e(u_1)e(u_2)-
\sh\,\pi\eta(u-i\tih)e(u_2)e(u_1)
&=&\sh(i\pi\eta\tih)
\left(e(u_1)^2+e(u_2)^2\right),\label{ee-univ}\\
\sh\,\pi\eta'(u-i\tih)f(u_1)f(u_2)-
\sh\,\pi\eta'(u+i\tih)f(u_2)f(u_1)
&=&-\sh(i\pi\eta'\tih)
\left(f(u_1)^2+f(u_2)^2\right),\label{ff-univ}\\
h(u_1)h(u_2)
{\sh\,\pi\eta(u+i\tih)\sh\,\pi\eta'(u-i\tih)
\over \sh\,\pi\eta'(u+i\tih)\sh\,\pi\eta(u-i\tih)}
&=&
h(u_2)h(u_1)\ .\label{hh}
\eea
\medskip

The proof is a direct substitution of the Gauss decomposition
of $L$-operators \r{GL-univ} into \r{RLL-univ}.

\subsection{The generating integrals  for $\Aelssect$}

Let $e^\pm(u)$, $f^\pm(u)$ and $h^\pm(u)$ be the following formal
 integrals of the symbols $\he_\la$, $\hf_\la$, $\hk_\la$ 
($u\in\CC$):
\bea
e^\pm(u)&=&\pm{\sin\,\pi\eta\tih\over\pi\eta}
\int_{-\infty}^\infty d\la\
\ee^{i\la u}\
{\he_\la \ee^{\mp c\tih\la/4}\over 1+\ee^{\pm\la/\eta}}\ ,
\label{Lapl-e}\\
f^\pm(u)&=&\pm{\sin\,\pi\eta'\tih\over\pi\eta'}
\int_{-\infty}^\infty d\la\  \ee^{i\la u}\
{\hf_\la \ee^{\pm c\tih\la/4}
\over 1+\ee^{\pm\la/\eta'}}\ ,
\label{Lapl-f}\\
h^\pm(u)
&=&-{\sin\,\pi\eta\tih\over2\pi\eta}
\int_{-\infty}^\infty d\la\  \ee^{i\la u}\
\hk_\la \ee^{\mp\la/2\eta''}\ ,
\label{Lapl-h}
\eea
Here
$$%\beq
\eta''={2\eta\eta'\over\eta+\eta'}\ .
$$%\label{eta''}\eeq
By the direct verification we can check that if the complex number
$u$ is inside the strip
$$%\beq
\Pi^+=\left\{-{1\over\eta}-{\tih c\over4}<\Im u<-{\tih c\over4}\right\}
$$%\label{strip+}\eeq
then elements $e^+(u)$, $f^+(u)$, $h^+(u)$ belong to $\Ael$. If the
complex number $u$ is inside the strip
$$%\beq
\Pi^-=\left\{{\tih c\over4}<\Im u<{\tih c\over4}+{1\over\eta}\right\}
$$%\label{strip-}\eeq
then the elements
\beq
e^-(u)=-e^+(u-i/\eta'')\ ,\quad
f^-(u)=-f^+(u-i/\eta'')\ ,\quad
h^-(u)=h^+(u-i/\eta'')\   \label{-+relat}
\eeq
also belong to $\Ael$. Thus we can treat the integrals
 $e^\pm(u)$, $f^\pm(u)$, $h^\pm(u)$
as generating functions of the elements of the algebra $\Ael$, analytical
in the strips $\Pi^{\pm}$. We  can state the following
\smallskip

\noindent
{\bf Proposition 4.} {\it
The generating functions $e(u)=e^\pm(u)$, $f(u)=f^\pm(u)$, $h(u)=h^\pm(u)$
satisfy the commutation relations}
\r{ef}--\r{hh} {\it if $\he_\la$, $\hf_\la$, $\hh_\la$ satisfy the 
relations} \r{he-hf}--\r{hh-hh}.
  
\smallskip

In order to prove this proposition we should use 
the Fourier transform calculations and fix in  \r{ef}--\r{hh}
either $\Im u_1<\Im u_2$ or $\Im u_1>\Im u_2$.

The relations \r{-+relat} and analiticity of generation functions in the 
 domains $\Pi^\pm$ allow one to make an analitical continuation of the 
 relations \r{ef}--\r{hh} including all possible combinations of 
the generating
 integrals. For instance, from (\ref{he}) we have also ($u=u_1-u_2$)
\bea
\sh\,\pi\eta(u+i\tih)h^\pm(u_1)e^\pm(u_2)&-&
\sh\,\pi\eta(u-i\tih)e^\pm(u_2)h^\pm(u_1)=\nn\\
&=&\sh(i\pi\eta\tih)
\{h^\pm(u_1),e^\pm(u_1)\}\ ,\nn\\
\sh\,\pi\eta(u+i\tih +i\tih c/2)h^+(u_1)e^-(u_2)&-&
\sh\,\pi\eta(u-i\tih +i\tih c/2)e^-(u_2)h^+(u_1)=\nn\\
&=&\sh(i\pi\eta\tih)
\{h^+(u_1),e^+(u_1)\}\ ,\nn\\
\sh\,\pi\eta(u+i\tih -i\tih c/2)h^-(u_1)e^+(u_2)&-&
\sh\,\pi\eta(u-i\tih -i\tih c/2)e^+(u_2)h^-(u_1)=\nn\\
&=&\sh(i\pi\eta\tih)
\{h^-(u_1),e^-(u_1)\}\ .\nn
\eea
Let now 
$$%\beq
R^-(u)=\tau^-(u)R(u), \quad \tau^-(u)=\th\left({\pi u\over2\tih}\right)
$$%\label{Runiv--}\eeq
and $e^\pm(u)$, $f^\pm(u)$, $h^\pm(u)$ be the Gauss coordinates of the 
 $L$-operators $L^\pm(u)$:
$$%\beq
L^\pm(u)
=\left(\begin{array}{cc} 1& f^\pm(u)\\0&1\end{array}\right)
\left(\begin{array}{cc}  (k^\pm(u+i\tih))^{-1}&0\\ 0&k^\pm(u)\end{array}\right)
\left(\begin{array}{cc} 1&0\\ e^\pm(u)&1\end{array}\right)\ .
$$%\label{GL-univ-pm}\eeq 
One can prove in an analogous manner that the described above mixed relations
 for generating functions $e^\pm(u)$, $f^\pm(u)$, $h^\pm(u)$ are equivalent to
 the following system of equations for the Gauss coordinates of 
the $L$-operators
($u=u_1-u_2$):
\begin{eqnarray}
R^+(u-ic\tih/2,\eta')L^+_1(u_1,\eta)L^-_2(u_2,\eta)&=&
L^-_2(u_2,\eta)L^+_1(u_1,\eta) R^+(u+ic\tih/2,\eta),\nn\\
R^\pm(u,\eta')L^\pm_1(u_1,\eta)L^\pm_2(u_2,\eta)&=&
L^\pm_2(u_2,\eta)L^\pm_1(u_1,\eta) R^\pm(u,\eta).
\label{RLL-univ-+-}
\eea
  These equations
 can be also obtained by means of the formal analytical
continuation  of the  relations \r{RLL-univ}.

The relation \r{-+relat} and the involution \r{auto} in the algebra $\Ael$
in terms of the $L$-operators
can be written as follows
\beq
L^+(u-i/\eta'')=
\sigma_z   L^-(u)    \sigma_z = \iota\sk{L^-(u)}
\label{L-oper-rel}
\eeq
and possibility to obtain  the commutation relations \r{RLL-univ-+-}
from \r{RLL-univ} by means of the analytical conti\-nua\-ti\-on follows
from
the quasi-periodicity property of the $R$-matrices $R^\pm(u,\eta)$:
$$%\beq
R^+(z-i/\eta)=
(\sigma_z\ot 1)   R^-(z)       (\sigma_z\ot 1)=
(1\ot \sigma_z)   R^-(z)       (1\ot\sigma_z)\ .
$$%\label{quasi-per}\eeq

Let us note that the
formal generators $\he_\la$, $\hf_\la$, $\hh_\la$
of the algebra $\Ael$ can be expressed through their generating
integrals using the inverse integral transform
\bea
\he_\la&=&\pm{\eta\ee^{\pm c\tih\la/4}\left(1+\ee^{\pm\la/\eta}\right)
\over 2\sin\,\pi\eta\tih}\int_{\Gamma^\pm}du\
\ee^{-i\la u} e^\pm(u)\ ,\nn\\
\hf_\la&=&\pm{\eta'\ee^{\mp c\tih\la/4}\left(1+\ee^{\pm\la/\eta'}\right)
\over 2\sin\,\pi\eta'\tih}\int_{\Gamma^\pm}du\
\ee^{-i\la u} f^\pm(u)\ ,\nn\\
\hh_\la&=&-{\eta\ee^{\pm \la/2\eta''}
\over \sin\,\pi\eta\tih}\int_{\Gamma^\pm}du\
\ee^{-i\la u} h^\pm(u)\ ,\label{inver}
\eea
where $\Gamma^\pm$ are contours which go from $-\infty$ to $+\infty$
inside the strips $\Pi^\pm$.

Using the 
relations (\ref{inver}) one can verify that the defining relations 
(\ref{c-ever})--(\ref{hh-hh}) 
for the algebra $\Ael$ are equivalent to the relations 
(\ref{ef})--(\ref{hh}) on generating functions of the algebra.
\medskip 

{\bf Remarks.}

{\bf 1.}
 The formal algebra generated by Gauss coordinates $e(u)$, $f(u)$ and $h(u)$
is not completely equivalent to the algebra of coefficients of $L(u)$ with
 the relations (\ref{RLL-univ})--(\ref{qdet=1}) since $h(u)$ is a quadratic
expression of $k_2(u)$. Naturally, one may consider 
the corresponding extension of $\Ael$, which looks a bit more complicated.
 Nevertheless, the algebra $\Ael$ is sufficient for the description of 
 representations in which we are interested in.

{\bf 2.} Let $\eta=1/\xi$. The matrix
\beq
S(\beta,\xi)=-(\sigma_z\otimes1)
R (\beta,1/\xi)
(1\otimes\sigma_z)
\label{trans}
\eeq
was obtained in \cite{ZZ} as an exact $S$-matrix of soliton-antisoliton
scattering in  Sine-Gordon model, where $\xi$ is related to
the coupling constant
of the model 
and we should set $\tih=\pi$ (we prefer 
to keep the free parameter $\tih$ for the convenience of taking the classical
limit \cite{KLPST}).
This $S$-matrix satisfies the conditions
of unitarity and crossing symmetry
$$%\beq\label{unitar}
S(\beta,\xi)S(-\beta,\xi)=1 \ ,
$$%\eeq
$$%\beq\label{cross}
(C\ot\id)\, S(\beta,\xi)\, (C\ot\id) =
\left(S (i\pi-\beta,\xi)\right)^{t_1}
$$%\eeq
with the charge conjugation matrix
$$C=\left(\begin{array}{cr} 0&1\\1&0 \end{array}\right).$$
Because of the relation \r{trans} the $R$-matrix \r{R-mat}
satisfies the same properties of unitarity and crossing
symmetry but with a different charge conjugation matrix
$$%\beq\label{charge}
\tilde C=\left(\begin{array}{cr} 0&1\\-1&0 \end{array}\right).
$$%\eeq
Since the square of this matrix is equal to $-1$
we have to use an unusual definition of the intertwining operators
\r{inteq} in order to have a possibility to identify them as
Zamolodchikov--Faddeev operators for the Sine-Gordon model. We will
discuss this  point in more details in the last section.

\subsection{The Yangian limit}

As follows from the definition of the elements of the algebra 
$\Ael$ \r{elements} and the generating functions \r{Lapl-e}--\r{Lapl-f}
each substrip of the strips $\Pi^\pm$ defines a subalgebra of the
algebra $\Ael$. In terms of the Fourier components these subalgebras
defined by  different asymptotics of the functions 
$\phi(\{\la_k\};\{\mu_i\};\{\nu_j\})$ at
$\la_k,\mu_i,\nu_j\to\pm\infty$ in \r{elements}.

Let us consider the substrips $\overline{\Pi}^\pm\subset \Pi^\pm$
$$
\overline{\Pi}^+=\left\{-{1\over2\eta}<\Im u<-{\tih c\over4}\right\}
,\quad
\overline{\Pi}^-=\left\{{\tih c\over4}<\Im u<{1\over2\eta}\right\}
$$
and restrict the generating functions  $e^\pm(u)$,  $f^\pm(u)$,  $h^\pm(u)$
onto these strips. Then in the limit $\eta\to0$ these generating
functions will be defined in the lower and the upper  half-planes and 
the relations \r{-+relat} drops out.
The
expressions via formal generators turn into the Laplace transform
(see Fig.~1.). For example,
$$
\left.e^\pm(u)\right|_{\eta=0}=\pm\tih\int_0^\infty d\la\
\ee^{\mp i\la (u\pm ic\tih\la/4)}\he_{\mp\la}\ ,
$$
and the defining relations \r{RLL-univ-+-} turn into the defining
relations of the central extended Yangian double \cite{K}.

\unitlength 1mm
\linethickness{0.4pt}
\begin{picture}(140.00,60.00)
\put(139.67,10.00){\vector(1,0){0.2}}
\put(40.00,10.00){\line(1,0){99.67}}
\bezier{120}(40.33,11.00)(70.00,11.00)(70.00,11.00)
\bezier{96}(70.00,11.00)(86.00,12.33)(90.00,19.00)
\bezier{96}(90.00,19.00)(94.67,26.33)(110.00,27.00)
\bezier{120}(110.00,27.00)(123.67,28.33)(140.00,28.00)
\put(90.00,10.00){\line(0,1){5.00}}
\put(90.00,17.67){\line(0,1){4.33}}
\put(90.00,26.33){\line(0,1){3.8}}
\put(90.00,30.00){\line(1,0){5.00}}
\put(100.00,30.00){\line(0,1){0.00}}
\put(100.00,30.00){\line(1,0){5.00}}
\put(110.00,30.00){\line(1,0){5.00}}
\put(120.00,30.00){\line(1,0){5.00}}
\put(130.00,30.00){\line(1,0){7.33}}
\put(139.67,40.00){\vector(1,0){0.2}}
\put(140.00,40.00){\line(-1,0){99.67}}
\bezier{120}(139.67,41.00)(110.00,41.00)(110.00,41.00)
\bezier{96}(110.00,41.00)(94.00,42.33)(90.00,49.00)
\bezier{96}(90.00,49.00)(85.33,56.33)(70.00,57.00)
\bezier{120}(70.00,57.00)(56.33,58.33)(40.00,58.00)
\put(90.00,40.00){\line(0,1){5.00}}
\put(90.00,47.67){\line(0,1){4.33}}
\put(90.00,56.33){\line(0,1){3.8}}
\put(90.00,60.00){\line(-1,0){3.00}}
\put(80.00,60.00){\line(0,1){0.00}}
\put(80.00,60.00){\line(-1,0){5.00}}
\put(70.00,60.00){\line(-1,0){5.00}}
\put(60.00,60.00){\line(-1,0){5.00}}
\put(50.00,60.00){\line(-1,0){7.33}}
\put(16.33,20.00){\makebox(0,0)[cc]{\mbox{${1\over1+\ee^{-\la/\eta}}$}}}
\put(16.33,50.00){\makebox(0,0)[cc]{\mbox{${1\over1+\ee^{\la/\eta}}$}}}
\put(16.33,5.00){\makebox(0,0)[cc]{Fig.~1.}}
\end{picture}

Although the commutation relation in terms of the generating
functions for the limiting algebra 
${{\cal A}_{\tih,0}(\widehat{\frak{sl}_2})}$
coincide with the commutation relations of the
central extended Yangian double $\widehat{DY(sl_2)}$
these two algebras should not be treated as isomorphic.
For instance, the algebra $\widehat{DY(sl_2)}$  has a discrete
set of generators and ${{\cal A}_{\tih,0}(\widehat{\frak{sl}_2})}$
 has a continious family of them. As a consequence, they have different 
 representation theories. It was also pointed out in \cite{Ko}. 
See details in \cite{KLP2}.

\subsection{Comultiplication structure}

As we already mentioned in the introduction, the algebra $\Ael$ is {\it not}
a Hopf algebra in the usual sense. Nevertheless we can assign the
Hopf algebra structure to the family of the algebras $\Ael$
parametrized by the parameter $\eta$.
Let us describe this Hopf structure.
In this subsection it is convenient to use instead of the
parameter $\eta$ its inverse
$$
\xi={1\over\eta}\ .
$$
Because of the relation \r{L-oper-rel} we can define the coproduct
only for one type of the operator, say, $L^+(u)$.

Consider the operation
\bea
\Delta\,c&=&c'+c''=c\otimes1+1\otimes c\ ,\nn\\
\Delta' L^+(u,\xi)&=&L^+(u-i\tih c''/4,\xi+\tih c'')
{\dot\otimes} L^+(u+i\tih c'/4,\xi)
\label{comul-L-univ}
\eea
or in components
\beq
\Delta L^+_{ij}(u,\xi)=
\sum_{k=1}^2
L^+_{kj}(u+i\tih c''/4,\xi) \otimes L^+_{ik}(u-i\tih c'/4,\xi+\tih c')
\label{com-L*cmp}
\eeq
which defines the coassociative map
$$%\beq
\Delta: \; \Aelx\rightarrow \Aelx \otimes
{\cal A}_{\tih,\xi+\tih c'}(\widehat{\frak{sl}_2})
$$%\label{map}\eeq
on the family of algebras $\Aelx$. The map $\Delta$ 
is a morphism of algebras, but it sends one algebra to a tensor product
 of two different algebras, which we do not identify. So we say that $\Ael$
 form (over the parameter $\eta$) a Hopf family of algebras.
Let us also note that because of the relation \r{L-oper-rel}
the comultiplication of the $L$-operator $L^-(u)$ is given by
$$
\Delta L^-_{ij}(u,\xi)=
\sum_{k=1}^2
L^-_{kj}(u-i\tih c''/4,\xi) \otimes L^-_{ik}(u+i\tih c'/4,\xi+\tih c')\ .
$$

In order
to save notations
below in this subsection and in the Appendix A we will understand
the operators $k(u,\xi)$, $e(u,\xi)$, $f(u,\xi)$ as operators
$k^+(u,\xi)$, $e^+(u,\xi)$, $f^+(u,\xi)$ and
will write $u',\xi'$ ($u'',\xi''$) in left (right)
components of the tensor product and understand them as
$u'=u+i\tih c''/4$, $\xi'=\xi=1/\eta$
($u''=u-i\tih c'/4$, $\xi''=\xi+\tih c'$).

The comultiplications of the operators $e(u,\xi)$, $f(u,\xi)$
and $h(u,\xi)$ are
\bea
\Delta e(u,\xi) &=&
e(u',\xi)\ot 1 +
\sum_{p=0} ^{\infty}(-1)^p
       \left(f(u'-i\tih,\xi')\right)^{p} h(u',\xi')\otimes
\left(e(u'',\xi'') \right)^{p+1}   \label{com-e-fu}\\
\Delta f(u,\xi)&=&
1\otimes f(u'',\xi'') +
\sum_{p=0} ^{\infty} (-1)^p
       \left(f(u',\xi')\right)^{p+1} \otimes
\tilde h(u'',\xi'')\left(e(u''-i\tih,\xi'')
\right)^p\ ,   \label{com-f-fu}\\
\Delta h(u,\xi)&=&\sum_{p=0}^\infty (-1)^p [p+1]_\eta
\left(f(u'-i\tih,\xi')\right)^{p} h(u',\xi') \ot
h(u'',\xi'')\left(e(u''-i\tih,\xi'') \right)^p
\label{comul-h}
\eea
where we define
$$%\beq
[p]_\eta={\sin\,\pi\eta\tih p\over \sin\,\pi\eta\tih}\ .
$$%\label{num-qua}\eeq
The proof of these formulas  is shifted to the Appendix A. Note that the 
 involution $\iota$ is compatible with coalgebraic structure: $\Delta \iota = 
(\iota\otimes \iota) \Delta$.

\setcounter{equation}{0}
\section{The Algebra $\Aelsect$ ($c=0$)}

Consider the formal algebra of the symbols $\he_\la$,
$\hf_\la$, $\hhh_\la$ and $S_0$ which satisfy the commutation
relations
\bea
{[}\he_\la,\hf_\mu{]}&=& \hhh_{\la+\mu}\ ,\label{he-hf0}\\
{[}S_0,\he_\mu{]}&=&
\sin(\pi\eta\tih)\,\tg(\pi\eta\tih)
\{\hhh_0,\he_\mu\}\ , \label{S-he0}\\
{[}\hhh_\la,\he_\mu{]}&=&
{\{S_0,\he_{\la+\mu}\} 
\over  \cos\,\pi\eta\tih} +\nn\\
&+&
{\tg\,\pi\eta\tih\over 2 \pi\eta}
\vpint d\tau\
\left[
\cth\left(\fract{\tau}{2\eta}\right)-\cth\left(\fract{\la+\tau}{2\eta}\right)
\right]
\{\hhh_{\la+\tau},\he_{\mu-\tau}\}\ ,
\label{hh-he0}\\
{[}S_0,\hf_\mu{]}&=& -
\sin(\pi\eta\tih)\,\tg(\pi\eta\tih)
\{\hhh_0,\hf_\mu\}\ , \label{S-hf0}\\
{[}\hhh_\la,\hf_\mu{]}&=& -
{\{S_0,\he_{\la+\mu}\} 
\over  \cos\,\pi\eta\tih} -\nn\\
&-&{\tg\,\pi\eta\tih\over 2 \pi\eta}
\vpint d\tau\
\left[
\cth\left(\fract{\tau}{2\eta}\right)-\cth\left(\fract{\la+\tau}{2\eta}
\right) \right]
\{\hhh_{\la+\tau},\hf_{\mu-\tau}\},
\label{hh-hf0}\\
{[}\he_\la,\he_\mu{]}&=&
{\tg\,\pi\eta\tih\over 2 \pi\eta}
\vpint d\tau\ \cth\left(\fract{\tau}{2\eta}\right)
\{\he_{\la+\tau},\he_{\mu-\tau}\}\ ,
\label{he-he0}\\
{[}\hf_\la,\hf_\mu{]}&=&
-{\tg\,\pi\eta\tih\over 2 \pi\eta}
\vpint d\tau\ \cth\left(\fract{\tau}{2\eta}\right)
\{\hf_{\la+\tau},\hf_{\mu-\tau}\}\ ,
\label{hf-hf0}\\
{[}S_0,\hhh_\la{]}&=&0\ ,\quad [\hhh_\la,\hhh_\mu]\ =\ 0\ ,
\label{hh-hh0}\\
1&=&S^2_0+\sin^2(\pi\eta\tih)\hhh_0^2\ .\label{addition}
\eea
We can assign to these commutation relations the analogous
ordering sense as we did it in the section 2.

Consider the free algebra $\overline{{\cal A}}_0$
formed by the formal integrals of the type
$$%\beq
\intt \prod_k d\la_k \prod_i d\mu_i \prod_j d\nu_j\
\phi_0(\{\la_k\};\{\mu_i\};\{\nu_j\})
P(\{\he_{\la_k}\};\{\hf_{\mu_i}\};\{\hhh_{\nu_j}\})
$$%\label{elements0}\eeq
where $\phi_0(\{\la_k\};\{\mu_i\};\{\nu_j\})$
 is the  $\CC$-number function of real
variables $\la_k$, $\mu_i$ and $\nu_j$ which
satisfy  the conditions
of analyticity:
\bea
&\phi_0(\{\la_k\};\{\mu_i\};\{\nu_j\})
&\ \ \mbox{is analytical in the strip}\ \ 
-\pi\eta<\Im(\la_k,\mu_i)<\pi\eta\ \ \forall\ \ \la_k,\mu_i\ ,\nn\\
&\phi_0(\{\la_k\};\{\mu_i\};\{\nu_j\})
&\ \ \mbox{is analytical in the strip}\ \ 
-2\pi\eta<\Im\nu_j<2\pi\eta\ \  \forall\ \ \nu_j\ ,\nn
\eea
except the points $\nu_j=0$ where this function has simple pole with
respect to all $\nu_j$.
The function
$\phi_0(\{\la_k\};\{\mu_i\};\{\nu_j\})$ has
the asymptotics when $\Re\la_k$, $\Re\mu_i$,
$\Re\nu_j\to\pm\infty$:
\bea
\phi_0(\{\la_k\};\{\mu_i\};\{\nu_j\})&<& C\ee^{-\alpha|\Re\la_k|}\ ,\nn\\
\phi_0(\{\la_k\};\{\mu_i\};\{\nu_j\})
&<& C\ee^{-\beta|\Re\mu_i|}\ ,\nn\\
\phi_0(\{\la_k\};\{\mu_i\};\{\nu_j\})
&<& C\ee^{-\gamma|\Re\nu_j|}\ ,\nn
\eea
for some real positive $\alpha$, $\beta$, $\gamma$.

The algebra $\Ael$ at $c=0$ is identified with $\overline{{\cal A}}_0$
factorized by the ideal generated by the commutation relations
\r{he-hf0}--\r{addition}.

Consider the formal integrals
\bea
e^+(u)&=&{\sin\,\pi\eta\tih\over\pi\eta}
\int_{-\infty}^\infty d\la\
\ee^{i\la u}\
{\he_\la \over 1+\ee^{\la/\eta}}\ ,
\label{Lapl-e0}\\
f^+(u)&=&{\sin\,\pi\eta\tih\over\pi\eta}
\int_{-\infty}^\infty d\la\  \ee^{i\la u}\
{\hf_\la
\over 1+\ee^{\la/\eta}}\ ,
\label{Lapl-f0}\\
h^+(u)
&=&S_0+{\sin\,\pi\eta\tih\over \pi\eta}
\int_{-\infty}^\infty d\la\  \ee^{i\la u}\
{\hhh_\la\over  1-\ee^{\la/\eta}}\ ,
\label{Lapl-h0}
\eea
as generating integrals of the elements of the algebra
$\Ael$ at $c=0$. We can prove that these generating functions are 
 analytical in the strip $\Pi^+=\{-1/\eta<\Im u<0\}$ and satisfy the
commutations relations \r{ef}--\r{hh}, where at $c=0$ we should
set $\eta'=\eta$. In particular, in this case the generating
functions $h^+(u)$ commute.

The different presentation of the generating functions $h^\pm(u)$ for
$c\neq 0$ (\ref{Lapl-h}) and for $c = 0$ (\ref{Lapl-h0}) follows from the
analysis of their asymptotical behaviour. Indeed, the relations 
(\ref{ef})--(\ref{hh})  imply that 
$$%\beq
e^+(u)\stackreb{\sim}{\Re u\to\pm\infty}
\ee^{-\pi\eta|u|},\quad
f^+(u)\stackreb{\sim}{\Re u\to\pm\infty}
\ee^{-\pi\eta|u|}\ ,
$$%\label{asympEF}\eeq
\beq
h^+(u)\stackreb{\sim}{\Re u\to\pm\infty}
\ee^{-\pi(\eta -\eta')|u|},
\label{asympHH}
\eeq
and constant, but different for $+\infty$ and $-\infty$ asymptotics of
$h^+(u)$ for $c = 0$:
$$%\beq
h^+(u)\stackreb{\sim}{\Re u\to\pm\infty}
 h^+(\pm\infty)\equiv h_\pm=S_0\mp i\sin(\pi\eta\tih)\hhh_0\ .
$$%\label{asympH}\eeq
Such asymptotics can be achieved by the following Cauchy kernel
presentations: 
\bea
e^+(u)&=&{\sh\,i\pi\eta\tih\over 2\pi}
\intt d\tilde  v\ { E(\tilde
v)\over\sh\,\pi\eta(\tilde  v-u)}\ ,
\label{2E+0}\\
f^+(u)&=&{\sh\,i\pi\eta\tih\over 2\pi} \intt d\tilde  v\ { F(\tilde
v)\over\sh\,\pi\eta(\tilde  v-u)}\ ,
\label{2F+0}\\
h^+(u)&=&S_0+ {\sh\,i\pi\eta\tih\over 2\pi}  \intt d\tilde
v\ H(\tilde  v)\,\cth\,\pi\eta(\tilde v-u)\ ,
\label{2H+0}
\eea
where $u \in \Pi^+=\{-1/\eta<\Im u<0\}$ and $\tilde v\in\RR$.
 Analogous  formulas take place for the ``$-$''-generating functions
 but with the spectral parameter $u$ in the strip
$\Pi^-=\{0<\Im u<1/\eta\}$. The presenations (\ref{2E+0})--(\ref{2H+0})
 are equivalent to the deformed Laplace
presentations 
(\ref{Lapl-e0})--(\ref{Lapl-h0}) if $E(\tilde{u})$, $F(\tilde{u})$
and $H(\tilde{u})$ are Fourier transforms of  $\he_\la$, $\hf_\la$ and 
$\hhh_\la$: 
$$
E(\tilde v) =
\int_{-\infty}^\infty d\la\ \ee^{i\la \tilde  v}  \he_\la,\quad
F(\tilde v) =
\int_{-\infty}^\infty d\la\ \ee^{i\la \tilde  v}  \hf_\la,\quad
H(\tilde v) =
\int_{-\infty}^\infty d\la\ \ee^{i\la \tilde  v}  \hhh_\la\ .
$$
For $c\neq 0$, due to \r{asympHH}, we use instead of \r{Lapl-h0} and 
\r{2H+0} the usual Fourier transform \r{Lapl-h}. Note that in the limit 
 $c\to 0$ we have the relation 
\beq
\hhh_\la=\hh_\la\ \sh\left(\fract{\la}{2\eta}\right)\ 
\label{new-old}
\eeq
 for all $\la$, so $\hhh_0$ is well defined whereas $\hh_0$ tends to infinity
 when $c\to 0$; to the contrary, $\hhh_0$ from \r{new-old} is zero for 
 $c\neq 0$ while $\hh_0\neq 0$ in this case.

The asymptotic generators $h_\pm$ of $c =0$ algebra $\Ael$
have the following commutation relations with the
generating functions $e^+(u)$ and $f^+(u)$:
$$%\beq
h_\pm e^+(u) h_\pm^{-1}=q^{\mp 2}e^+(u),\quad
h_\pm f^+(u) h_\pm^{-1}=q^{\pm 2}f^+(u),\quad
q=\ee^{i\pi\eta\tih}\ .
$$%\label{q-comm}\eeq
and are primitive group-like elements: $\Delta h_\pm=h_\pm\otimes h_\pm$.
Thus their product 
 $$h_+ h_- = S^2_0+\sin^2(\pi\eta\tih)\hhh_0^2$$
 is cental and group-like primitive. Due to this we can put it to be equal to
1. This kills unnecessary representations of level 0.

\subsection{Evaluation homomorphism}

Let $e$, $f$ and $h$ be the generators of the algebra 
$U_{i\pi\eta\tih}(sl_2)$:
$$ [h,e]=2e,\qquad [h,f]=-2f,\qquad [e,f]=[h]_\eta= 
{\sin\,\pi\eta\tih h\over \sin\,\pi\eta\tih}\ .$$
The following proposition presents two descriptions of the evaluation
homomorphism of $c=0$ algebra $\Ael$ onto
 $U_{i\pi\eta\tih}(sl_2)$.
\smallskip

\noindent
{\bf Proposition 5.} {\it The algebra $\Ael$ at $c=0$
has the following  evaluation homomorphism $\Ev_z$ onto $U_q({sl}_2)$},
 $z\in \CC$
\bea
\Ev_z \sk{S_0}&=&\cos(\pi\eta\tih h)\ ,\nn\\
\Ev_z \sk{\he_\la}&=&\ee^{-i\la z}\, \ee^{-\tih\la(h-1)/2}\,  e =
\ee^{-i\la z}\, e\, \ee^{-\tih\la(h+1)/2}\ ,\nn\\
\Ev_z \sk{\hf_\la}&=&\ee^{-i\la z}\, \ee^{-\tih\la(h+1)/2}\,  f =
\ee^{-i\la z}\, f\, \ee^{-\tih\la(h-1)/2}\ ,\nn\\
\Ev_z \sk{\hhh_\la}&=&\ee^{-i\la z}\, \ee^{-\tih\la(h-1)/2}\, ef
-\ee^{-i\la z}\, \ee^{-\tih\la(h+1)/2}\, fe\ \nn
\eea
{\it or, equivalently,} $(u\in\Pi^+)$
\bea
\Ev_z \sk{e^+(u)}&=&-
{\sh\,i \pi\eta \tih \over \sh\,\pi\eta (u-z+i\tih(h-1)/2))  }\, e =
- e\, {\sh\,i \pi\eta \tih \over \sh\,\pi\eta (u-z+i\tih(h+1)/2))  }\ ,
\nn\\
\Ev_z \sk{f^+(u)}&=&-
{\sh\,i \pi\eta \tih \over \sh\,\pi\eta (u-z+i\tih(h+1)/2))  }\, f =
- f\, {\sh\,i \pi\eta \tih \over \sh\,\pi\eta (u-z+i\tih(h-1)/2))  }\ ,
\nn\\
\Ev_z \sk{h^+(u)}&=&
\cos(\pi\eta\tih h) - \sh\,i\pi\eta\tih\left[
\cth\,\pi\eta(u-z+i\tih(h-1)/2)\, ef\right.\nn\\
&-&\left.\cth\,\pi\eta(u-z+i\tih(h+1)/2)\, fe
\right].\nn
%\label{eval-h-n}
\eea
\smallskip

Let $V_n$ be $(n+1)$-dimensional $U_{i\pi\eta\tih}(sl_2)$-module with a
 basis $v_k$, 
$k=0,1,\ldots,n$ where  the operators $h$, $e$ and $f$ 
 act  
according to the rules
$$%\beq
h\,v_k=(2k-n)\,v_k,\quad e\,v_k=[k]_\eta\,v_{k-1},\quad
f\,v_k=[n-k]_\eta\,v_{k+1}\ .
$$%\label{finite-n}\eeq
Due to the Proposition 5 we have an action $\pi_n(z)$ of the algebra $\Ael$ in 
 the space 
$V_{n,z} = V_n $.
Note that the action of $h^+(u)$ can be simplified in this case as
\beq
\pi_n(z)\left(h^+(u)\right)\ =\
{
\sh\,\pi\eta (u-z-i\tih(n+1)/2))
\sh\,\pi\eta (u-z+i\tih(n+1)/2))
\over
\sh\,\pi\eta (u-z+i\tih(h+1)/2))
\sh\,\pi\eta (u-z+i\tih(h-1)/2))
}\ .
\label{conven}
\eeq

The simplest two-dimensional evaluation representation $\pi_1(z)$
 of the
algebra $\Ael$ on the space 
$V_z =V_{1,z}$ (we have identified $v_{+,z}=v_{0}$
and $v_{-,z}=v_{1}$) is
\bea
e^+(u)v_{+,z}&=&0,\quad f^+(u)v_{-,z}\, =\, 0\ , \label{eval-0}\\
e^+(u)v_{-,z}&=& -
{\sh\,i \pi\eta \tih \over \sh\,\pi\eta (u-z)  }\ 
 v_{+,z}\ ,\quad
f^+(u)v_{+,z}\ =\ -
{\sh\,i \pi\eta \tih \over \sh\,\pi\eta (u-z)  }\ 
 v_{-,z}\ ,\label{eval-e-f1}\\
h^+(u)v_{\pm,z}&=&
\sh\,i\pi\eta\tih
\left[\cth\,i\pi\eta\tih\mp\cth\,\pi\eta(u-z)\right]v_{\pm,z}
= {\sh\,\pi\eta(u-z\mp i\tih) \over \sh\,\pi\eta (u-z)  }
 v_{\pm,z}
\ .\label{eval-h1}
\eea
The action of ``$-$''-generating functions on the space $V_z$ is
given by the same formulas \r{eval-0}--\r{eval-h1} but with $u\in\Pi^-$. 
In $L$-operator's terms, the representation described by
\r{eval-0}--\r{eval-h1} is equivalent to the standard one:
$$
\pi_1(z)L^\pm(u)=R^\pm(u-z,\eta).
$$

For $c=0$ in addition to the evaluation homomorphism we have, analogously
to the case of $U_q(\widehat{\frak{sl}_2})$, an embedding of a subalgebra,
 isomorphic to $U_q(sl_2)$. Here $q = \ee^{i\pi\eta\tih}$, $|q|= 1$.
This subalgebra is generated by the elements $S_0$, $\he_0$, $\hf_0$  
and $\hhh_0$ and is given in the form of Sklyanin degenerated algebra
 \cite{Sk}:
\bea
{[}\he_0,\hf_0{]}&=& \hhh_0\ ,\quad
{[}S_0,\hhh_0{]}\ =\ 0\ ,\quad
S_0^2+\sin^2(\pi\eta\tih)\hhh^2_0\ =\ 1\ ,
\nn\\
{[}S_0,\he_0{]}&=&
\sin(\pi\eta\tih)\,\tg(\pi\eta\tih)
\{\hhh_0,\he_0\}\ , \quad
{[}\hhh_0,\he_0{]}\ =\
{\tg\,\pi\eta\tih\over \sin\, \pi\eta\tih}
\{S_0,\he_0\}\ , \nn\\
{[}S_0,\hf_0{]}&=& -\sin(\pi\eta\tih)\,
\tg(\pi\eta\tih)
\{\hhh_0,\hf_0\}\ , \quad
{[}\hhh_0,\hf_0{]}\ =\ -
{\tg\,\pi\eta\tih\over \sin \pi\eta\tih}
\{S_0,\hf_0\}\ . \nn
\eea
But, unlikely to quantum affine case, this subalgebra is
  {\it not} a Hopf subalgbera. Moreover, this embedding is destroied
when $c \neq 0$. For instance, the genertors $\he_0$ and $\hf_0$
 commute for $c \neq 0$.

 In the rational limit $\eta\to0$  this
finite-dimensional subalgebra becomes $\frak{sl}_2$ subalgebra
of the Yangian double 
${{\cal A}_{\tih,0}(\widehat{\frak{sl}_2})}$.

\setcounter{equation}{0}
\section{Current Realization of $\Aelsect$}

In this section we would like to give another realization
of the algebra $\Ael$ which is an analog of the current realization
of the affine Lie algebras.
The necessity of this realization follows from the construction
of  infinite-dimensional representations of the algebra $\Ael$
at $c\neq0$ in terms of free fields.
  The  generating functions $e^\pm(u)$ and $f^\pm(u)$ cannot be realized
directly in terms of  free fields, but only some combinations of them, called
total currents, have a free field realization.

Let us define generating functions (total currents) $E(u)$ and   $F(u)$ 
 as
 formal Fourier transforms of the symbols $\he_\la$ and  $\hf_\la$ :
$$%\beq
E(u)=\intt d\la\ \ee^{i\la u} \he_\la,\quad
F(u)=\intt d\la\ \ee^{i\la u} \hf_\la,\quad u\in\CC
$$%\label{Fourier}\eeq
and put 
$$%\beq
H^\pm(u)= -{\tih\over 2}
\int_{-\infty}^\infty d\la\  \ee^{i\la u}\
\hk_\la \ee^{\mp\la/2\eta''},\qquad 
h^\pm(u)={\sin\,\pi\eta\tih\over\pi\eta\tih}\ H^\pm(u)\  .
$$%\label{H+-}\eeq
 We prove in this section that:

{\bf \it (i)} The currents 
 $E(u)$,  $F(u)$  and $H^\pm(u)$  satisfy the relations \r{h+h-c}--\r{efc}:
\bea
H^+(u) H^-(v)&=&
{
\sh\,\pi\eta(u-v-i\tih(1-c/2))
\sh\,\pi\eta'(u-v+i\tih(1-c/2))
\over
\sh\,\pi\eta(u-v+i\tih(1+c/2))
\sh\,\pi\eta'(u-v-i\tih(1+c/2))
}\ 
H^-(v) H^+(u)\ ,\label{h+h-c}\\
H^\pm(u) H^\pm(v)&=&
{
\sh\,\pi\eta(u-v-i\tih)
\sh\,\pi\eta'(u-v+i\tih)
\over
\sh\,\pi\eta(u-v+i\tih)
\sh\,\pi\eta'(u-v-i\tih)
}\
H^\pm(v) H^\pm(u)\ ,\label{hhc}\\
H^\pm(u)E(v)&=&
{
\sh\,\pi\eta(u-v-i\tih(1\mp c/4))
\over
\sh\,\pi\eta(u-v+i\tih(1\pm c/4))
}\
E(v)H^\pm(u)\ , \label{hec}\\
H^\pm(u)F(v)&=&
{
\sh\,\pi\eta'(u-v+i\tih(1\mp c/4))
\over
\sh\,\pi\eta'(u-v-i\tih(1\pm c/4))
}\
F(v) H^\pm(u)\ , \label{hfc}\\
E(u)E(v)&=&
{
\sh\,\pi\eta(u-v-i\tih)
\over
\sh\,\pi\eta(u-v+i\tih)
}\  E(v)E(u)\ ,\label{eec}\\
F(u)F(v)&=&
{
\sh\,\pi\eta'(u-v+i\tih)
\over
\sh\,\pi\eta'(u-v-i\tih)
}\ F(v)F(u)\ ,\label{ffc}\\
{[}E(u),F(v){]}&=& {2\pi\over\tih}\left[
\delta\left(u-v-\fract{ic\tih}{2}\right)
H^+\left(u-\fract{ic\tih}{4}\right)-
\delta\left(u-v+\fract{ic\tih}{2}\right)
H^-\left(v-\fract{ic\tih}{4}\right)
\right]  ,            \label{efc}
\eea
where the $\delta$-function is defined as
$$%\beq
2\pi \delta(u)=\lim{\epsilon\to0}{1\over i}
\left[ {1\over u-i\ep}-{1\over u+i\ep}\right]=
\intt d\la\ \ee^{i\la u},\quad u\in\RR\ .
$$%\label{delta}\eeq

{\bf\it (ii)} The generating functions $e^\pm(u)$ and $f^\pm(u)$
of the algebra $\Ael$ can be obtained from 
 $E(u)$ and $F(u)$ using the formulas:
\bea
e^\pm(u)&=&\sin\,\pi\eta\tih\int_{C} {dv\over2\pi i}\
 { E(v)\over\sh\,\pi\eta
(u-v\pm ic\tih/4)}\ ,\label{2e}\\
f^\pm(u)&=&\sin\,\pi\eta'\tih \int_{C'}
{dv\over2\pi i}\  { F(v)\over\sh\,\pi\eta'
(u-v\mp ic\tih/4)}\ ,\label{2f}
\eea
where the contour $C'$ goes from $-\infty$ to $+\infty$, 
the points $u+ic\tih/4+ik/\eta'$ ($k\geq0$) are above the contour and
the points $u-ic\tih/4-ik/\eta'$ ($k\geq0$) are below the contour.
The contour $C$ also goes from $-\infty$ to $+\infty$ but
the points $u-ic\tih/4+ik/\eta$ ($k\geq0$) are above the contour and
the points $u+ic\tih/4-ik/\eta$ ($k\geq0$) are below the contour.
In the considered case when $\tih>0$ the form of the contours
is shown on the Fig.~2.

\unitlength 1mm
\linethickness{0.4pt}
\begin{picture}(130.67,38.67)
\put(130.33,20.00){\vector(1,0){0.2}}
\put(10.33,20.00){\line(1,0){120.00}}
\bezier{148}(10.33,18.67)(42.67,18.67)(47.33,17.33)
\bezier{40}(47.33,17.33)(52.33,15.67)(56.33,12.33)
\bezier{44}(56.33,12.33)(59.67,9.33)(66.00,7.67)
\bezier{68}(66.00,7.67)(77.00,5.33)(83.00,7.00)
\bezier{84}(83.00,7.00)(92.67,12.33)(84.00,16.33)
\bezier{92}(84.00,16.33)(73.67,20.33)(62.00,24.00)
\bezier{80}(62.00,24.00)(53.33,28.67)(62.00,33.00)
\bezier{92}(62.00,33.00)(73.67,37.33)(83.00,33.33)
\bezier{48}(83.00,33.33)(91.00,30.00)(91.33,26.67)
\bezier{32}(91.33,26.67)(93.00,22.33)(96.33,22.00)
\bezier{76}(96.33,22.00)(105.67,21.00)(115.33,21.00)
\put(130.33,21.00){\vector(1,0){0.2}}
\put(115.33,21.00){\line(1,0){15.00}}
\put(74.33,30.00){\makebox(0,0)[cc]{$\bullet$}}
\put(74.33,11.33){\makebox(0,0)[cc]{$\bullet$}}
\put(60.67,16.00){\makebox(0,0)[lc]{$u-ic\tih/4$}}
\put(72.00,24.67){\makebox(0,0)[lc]{$u+ic\tih/4$}}
\put(130.67,25.67){\makebox(0,0)[cc]{$+\infty$}}
\put(10.33,26.00){\makebox(0,0)[cc]{$-\infty$}}
\put(10.00,9.33){\makebox(0,0)[cc]{Fig.~2.}}
\put(114.33,14.67){\makebox(0,0)[cc]{$C'$}}
\put(45.00,30.67){\makebox(0,0)[cc]{$C$}}
\end{picture}

{\bf\it (iii)}
The total currents $E(u)$ and $F(u)$ can be expressed through $e^\pm(u)$
 and $f^\pm(u)$ by means of the Ding-Frenkel relations \cite{FD}:
\bea
e^+\sk{u-\fract{ic\tih}{4}}-e^-\sk{u+\fract{ic\tih}{4}}&=&
{\sin\, \pi\eta\tih\over\pi\eta} E(u)\ ,\label{FDE}\\
f^+\sk{u+\fract{ic\tih}{4}}-f^-\sk{u-\fract{ic\tih}{4}}&=&
{\sin\, \pi\eta'\tih\over\pi\eta'} F(u)\ .\label{FDF}
\eea

 Let us start from {\it (ii)}. One can see that in every particular case of
 $e^+(u)$,  $e^-(u)$, $f^+(u)$ and $f^-(u)$ the contours of 
integration could be
 deformed to a straight line being the boundary of one of the strips $\Pi^\pm$.
 Then the relations \r{Lapl-e}--\r{Lapl-f} are equivalent to \r{2e}--\r{2f}
 via Fourier transform. Moreover, the relations \r{2e}--\r{2f} say that 
 $E(u)$ and $F(u)$ are the differences of boundary values of analitical 
  functions $e^\pm(u)$ and $f^\pm(u)$ for the Riemann problem on a strip.
Let us demonstrate this for $e^+(u)$.

When  the spectral parameter $u\in\Pi^+$  tends to the upper and
the lower boundaries
of the strip $\Pi^+$  we can obtain from \r{Lapl-e}:
\bea
\lim{\eps\to+0} (
e^+(\tilde u-ic\tih/4-i\eps)  + e^+(\tilde u-i/\eta-ic\tih/4+i\eps))
&=&{\sh\,i\pi\eta\tih\over i\pi\eta}
E(\tilde u),\quad \tilde u\in\RR\ ,
\nn\\
\lim{\eps\to+0} (
e^+(\tilde u-ic\tih/4-i\eps)  -
e^+(\tilde u-i/\eta-ic\tih/4+i\eps))
&=&  {\sh\,i\pi\eta\tih\over \pi}
\ \ \vpint d\tilde v\ { E(\tilde v)\over\sh\,\pi\eta
(\tilde v-\tilde u)}\ . \label{sohot2}
\eea
These relations are Sokhotsky-Plemely's formulas
associated with the Riemann problem
on the strips of the width $1/\eta$. Summing  the formulas \r{sohot2}
and using the analytical continuation with respect to the spectral
parameter $\tilde u$
we obtain \r{2e} where the contour of integration goes from $-\infty$
to $+\infty$ in such a way that point $u+i/\eta+ic\tih/4$ is above the
contour and the point $u+ic\tih/4$ is below. The same arguments
applied to the generating function $e^-(u)$ lead to \r{2e} but
with contour going between points $u-ic\tih/4$ and
$u-i/\eta-ic\tih/4$. 
%The loop-like contour $C$ shown on the Fig.~2
%can be explained by the requirement to have the same contour in the
%formulas \r{2e} for the generating functions  $e^+(u)$ and $e^-(u)$.
Repeating this consideration
for the generating functions $f^+(u)$ and $f^-(u)$ we obtain
%in 
\r{2f}. 
%the line-like contour  $C'$.

The relations \r{FDE}--\r{FDF} follow from \r{sohot2}. 
The commutation relations \r{h+h-c}--\r{efc} are direct corollaries of 
 \r{ef}--\r{hh} and \r{FDE}--\r{FDF}.
\medskip

{\bf Remark.} Note that, as usual for affine algebras, the relations 
 \r{FDE}--\r{FDF} should be understood in a sense of analytical continuation.
 For instance, if the argument of $e^+(u)$ is inside of its domain of 
 analiticity, the argument of $e^-(u)$ does not. It means, in particular,
 that the total currents $E(u)$ and $F(u)$ belong not the algebra $\Ael$ but
 rather to some its analytical extension. Nevertheless, they act on 
 highest weight representations, and the precise definition of the category of 
 highest weight representations should be equivalent to the description 
 of the proper analitical extension of the algebra $\Ael$.

\setcounter{equation}{0}
\section{Representations of the Algebra $\Aelsect$ at Level 1}

\subsection{Representation of the commutation relations by a free field}

This section is devoted to the construction of an infinite dimensional
representations of the algebra $\Ael$. For simplicity we will
consider the representation at level 1 ($c=1$) since in this case
only one free field is sufficient instead of three free fields for the
general case. The generalization to an arbitrary level could be done using
ideas developed in \cite{Mat,AOS}.
Here and till the end of the paper we
will understand everywhere $\eta'$ equal
$\eta/(1+\eta\tih)$.

Define bosons $a_\la$, $\la\in\RR$ which satisfy
the commutation relations \cite{L,JKM}:
\beq
[a_\la,a_\mu]={1\over\tih^2}\,{\sh(\tih\la)\sh(\tih\la/2)\over\la}\,
{\sh(\la/2\eta)\over\sh(\la/2\eta')}\, \delta(\la+\mu)=
\alpha(\la)\delta(\la+\mu)\ .
\label{bosons}
\eeq
Introduce also  bosons $a'_\la$ related to the initial ones as
follows:
$$%\beq
a'_\la={\sh(\la/2\eta')\over\sh(\la/2\eta)}a_\la\ .
$$%\label{bosons'}\eeq
Consider the generating functions
\bea
E(u)&=&\ee^{\gamma}\  
{:}\exp\left(\tih \intt d\la\ \ee^{i\la u}{a'_\la\over
\sh(\tih\la/2)}\right){:}\ ,\label{Ebos}\\
F(u)&=&\ee^{\gamma}\ {:}\exp\left(-\tih \intt d\la\ \ee^{i\la u}{a_\la\over
\sh(\tih\la/2)}\right){:}\ ,\label{Fbos}\\
H^\pm(u)&=&\ee^{-2\gamma}
\ {:}\ E\left(u\pm\fract{i\tih}{4}\right)
F\left(u\mp\fract{i\tih}{4}\right){:} =\nn\\
&=&  {:}\exp\left(\mp2\tih \intt d\la\ \ee^{i\la u}
{a_\la\ee^{\mp\tih\la/4}\over
1-\ee^{\pm\la/\eta}}\right){:}\ ,
\label{Hbos}
\eea
where $\gamma$ is the Euler constant and
the product of these generating functions is defined according
to \cite{JKM}
\bea
&{:}\exp\left(\intt d\la\ g_1(\la)\,a_\la\right){:}\ \cdot\
{:}\exp\left(\intt d\mu\ g_2(\mu)\,a_\mu\right){:} = \nn\\
&\quad =\exp\left(\int_{\tilde C}{d\la\,\ln(-\la)\over2\pi i}\
\alpha(\la)g_1(\la)g_2(-\la) \right)
{:}\exp\left(\intt d\la\ (g_1(\la)+g_2(\la))\,a_\la\right){:}\ .
\label{normal}
\eea
The contour $\tilde C$ is shown on the Fig.~3.
\bigskip

\unitlength 1.00mm
\linethickness{0.4pt}
\begin{picture}(121.00,20.00)
\put(17.00,15.00){\makebox(0,0)[cc]{0}}
\put(20.00,15.00){\makebox(0,0)[cc]{$\bullet$}}
\put(132.00,15.00){\makebox(0,0)[cc]{$+\infty$}}
\put(20.00,15.00){\line(1,0){100.33}}
\put(40.00,10.00){\line(1,0){80.33}}
\put(120.00,20.00){\line(-1,0){100.00}}
\put(30.00,5.00){\makebox(0,0)[cc]{Fig.~3.}}
\put(121.00,10.00){\vector(1,0){0.2}}
\put(100.00,10.00){\line(1,0){21.00}}
\put(20.00,10.00){\line(1,0){22.00}}
\put(20.67,10.00){\line(1,0){22.00}}
%\put(42.67,10.00){\line(0,1){0.00}}
\put(20.17,15.00){\oval(15.00,10.00)[l]}
\end{picture}
\smallskip

We have the following

\noindent
{\bf Proposition 6.} {\it The generating functions} \r{Ebos}--\r{Hbos}
{\it satisfy the commutation relations} \r{h+h-c}--\r{efc}.
\smallskip

The proof is based on the normal ordering relations gathered in the
Appendix B. We will show in the next subsection that the $\zeta$-function
regularization used in \cite{JKM} to define \r{normal} can be 
included into the definition  of the Fock space.

\subsection{Representation in a Fock space and twisted intertwining operators}

The goal of this subsection is to interpret Zamolodchikov-Faddeev
operators \cite{L,JKM} following the ideology \cite{DFJMN,JMbook}
as twisted intertwining operators for an 
infinite-dimensional representation of $\Ael$.
 For a description of this infinite-dimensional representation we need 
 a definition of a Fock space generated by continuous family of free bosons
 together with a construction of vertex operators.
We do it below in a slightly more general setting.
\bigskip

Let $a(\la)$ be a meromorphic function, regular for $\la\in \RR$ and satisfying
 the following conditions:
$$a(\la)= -a(-\la)\ ,$$
$$ 
a(\la) \sim a_0 \la,\quad \la \to 0, \qquad a(\la)\sim 
\ee^{a' |\la|},  \quad \la \to \pm\infty .
$$

Let $a_\la$, $\la\in\RR$, $\la\neq 0$ be free bosons which satisfy the 
 commutation relations
$$%\beq
[a_\la, a_\mu]=a(\la)\delta(\la+\mu) \ .
$$%\label{boson}\eeq

 We define a  (right) Fock space $\H_{a(\la )}$ as follows. 
$\H_{a(\la )}$ is generated as a vector
space 
by the expressions 
$$
\mint f_n(\la_n) a_{\la_n} d\la_n\ldots \mint f_1(\la_1) a_{\la_1} d\la_1\
\rvac\  , $$
 where the functions $f_i(\la)$ satisfy the condition
$$%\begin{equation}
f_i(\la)< C\, \ee^{(a'/2+\epsilon)\la}, \qquad \la \to -\infty\ ,
$$%\end{equation}
for some $\epsilon>0$ and 
$f_i(\la)$ are analytical functions in a neighbourhood of
 ${\RR_+}$ except $\la=0$, where they have a simple pole.

The left Fock space $\H^*_{a(\la)}$ is generated by the expressions
$$ 
\lvac \nint g_1(\la_1) a_{\la_1} d\la_1\ldots \nint g_n(\la_n) a_{\la_n} 
d\la_n\ ,
$$
where the functions $g_i(\la)$ satisfy the conditions
$$%\begin{equation}
g_i(\la)< C\, \ee^{-(a'/2+\epsilon)\la}, \qquad \la \to +\infty\ ,
$$%\end{equation}
for some $\epsilon>0$ 
and $g_i(\la)$ are analytical functions in a neighbourhood
of
 ${\RR_-}$ except $\la=0$, where they also have a simple pole.

The pairing $(,):$ $\H_{a(\la )}^*\otimes \H_{a(\la )} \to \CC$ is 
uniquely defined by
the
 following prescriptions:
\bea
 &(i)&\quad  (\langle\mbox{vac}|, \rvac) =1\ ,\nn\\
 &(ii)& \quad (\lvac \nint d\la \ g(\la)a_\la \ , 
\mint d\mu\  f(\mu)a_\mu\   \rvac) =
\int_{\tilde C} {d\la\,\ln(-\la)\over 2\pi i}
  g(\la)f(-\la)a(\la)\ ,\nn\\
&(iii)& \quad \mbox{the Wick theorem}.\nn
\eea

Let the vacuums $\lvac$ and $\rvac$ satisfy 
 the conditions
$$
a_\la\rvac =0,\quad \la>0,\qquad \lvac a_\la =0, \quad \la<0\ ,
$$
 and $f(\la)$ be  a function analytical in some neighbourhood of the real line
 with possible simple pole at $\la=0$ and which has the following asymptotical
 behaviour:
$$
\qquad f(\la)< e^{-(a'/2+\epsilon)|t|}
, \qquad \la \to \pm\infty
$$
for some $\epsilon>0$. 
Then, by definition, the operator
$$F= {:}\exp \sk{\int_{-\infty}^{+\infty}d\la\ f(\la)a_\la} {:}$$
acts on the right Fock space $\H_{a(\la )}$ as follows.
$F=F_-F_+$, where 
$$F_-=\exp \sk{\int_{-\infty}^{0}d\la\ f(\la)a_\la} \quad 
\mbox{and}\quad 
F_+=
\lim{\epsilon\to +0} \ee^{
\epsilon \ln \epsilon f(\epsilon)a_{\epsilon}}
\exp\sk{\int_{\epsilon}^{\infty}d\la\ f(\la)a_\la}\ .
$$
The action of operator $F$ on the left Fock space $\H_{a(\la )}^*$ is
 defined via another decomposition: $F=\tilde{F}_-\tilde{F}_+$, where 
$$\tilde{F}_+=\exp \sk{\nint d\la\  f(\la)a_\la}\quad
\mbox{and}\quad  
\tilde{F}_-=
\lim{\epsilon\to +0} \ee^{
\epsilon \ln\epsilon f(-\epsilon)a_{-\epsilon}}
\exp\sk{\int^{-\epsilon}_{-\infty}d\la\ f(\la)a_\la}.
$$
 These definitions imply the following statement:
\smallskip

\noindent
{\bf Proposition 7}.

$(i)$  {\it The defined above actions of the operator 
$$
F= :\exp \sk{\int_{-\infty}^{+\infty}d\la\ f(\la)a_\la}{:}
$$
on the Fock spaces $\H$ and $\H^*$ are adjoint};
 
$(ii)$ 
{\it The product of the normally ordered operators satisfy the property}
 \r{normal}.
\bigskip

Returning to level one representation of $\Ael$ we choose $\H = \H_{a(\la )}$
 for $a(\la)$ defined in \r{bosons}:
$$a(\la)={1\over\tih^2}\,{\sh(\tih\la)\sh(\tih\la/2)\over\la}\,
{\sh(\la/2\eta)\over\sh(\la/2\eta')}\ .$$
From the definition of the Fock space $\H$ and from the proposition 6 we have
 immediately the construction of a representation of $\Ael$:
\smallskip

\noindent
{\bf Proposition 8}. {\it The relations} \r{Ebos}--\r{Hbos} 
{\it define a highest
 weight level $1$ representation of the algebra $\Ael$ in the Fock space $\H$.}
\smallskip

The highest weight property means that  
$$%\beq
\he_\la\vac=0,\quad \hf_\la\vac=0,\quad \la>0\quad\mbox{and}\quad
\lvac\he_\la=0,\quad \lvac\hf_\la=0,\quad \la<0\ .
$$%\label{hwr}\eeq

Let us define the twisted intertwining operators
\begin{eqnarray}
\Phi(z)&:& \H\to \H\otimes V_{z+i\tih/2}\ ,\quad
\Phi^{*}(z)\ :\ \H\ot V_{z+i\tih/2} \to \H\ ,\nn%\label{typeI}
\\
\Psi^{*}(z)&:&V_{z+i\tih/2}\otimes \H\to \H\ ,\quad
\Psi(z)\ :\  \H\to  V_{z+i\tih/2}\ot \H\ , \nn%\label{typeII}
\end{eqnarray}
as those which commute with the action of the algebra $\Ael$ up to
the involution \r{auto}
\begin{eqnarray}
\Phi(z) \iota(x) &=& \Delta(x) \Phi(z)\ ,\quad
\Phi^{*}(z)\Delta(x)\ =\  \iota(x)\Phi^{*}(z)\ , \nn\\
\Psi^{*}(\b)\Delta(x)&=& \iota(x)\Psi^{*}(z)\ ,\quad
\Psi(z) \iota(x) \ =\  \Delta(x) \Psi(z)\ ,\quad
                          \label{inteq}
\end{eqnarray}
for arbitrary $x\in\Ael$.  In the definition of the intertwining
operators $V_{z}$ denotes the two-dimensional
evaluation module
$$V_{z}=V\otimes \CC[[\ee^{i\la z}]]\quad\mbox{and}\quad
V=\CC\, v_+\oplus\CC\, v_-\ , \quad \la\in\RR,\quad z\in\CC\ .
$$

The components of the intertwining operators are defined as follows:
\bea
\Phi(z) v  &=& \Phi_+(z)v\otimes v_+ + \Phi_-(\z)v\otimes v_-\ ,\quad
\Phi^{*}(z)(v\ot v_\pm)\ =\  \Phi^{*}_\pm(z)v \ , \nn\\
\Psi^{*}(z)(v_\pm\otimes v)&=& \Psi^{*}_\pm(z)v \ ,\quad
\Psi(z) v  \ =\  v_+\ot \Psi_+(z)v + v_-\ot \Psi_-(\z)v\ , \nn
\eea
where $v\in \H$ and one should understand the components
$\Phi_\ep(\z)$, $\Psi^{*}_\ep(\b)$,
$\ep=\pm$ as generating functions, for example:
$$%\beq
\Phi(\z)v=\sum_{\ep=\pm} \intt d\la\
\Phi_{\ep,\la}v\otimes v_{\ep}\ee^{i\la (z+i\tih/2)}\ .
$$%\label{genfun}\eeq
%The integral transform in \r{genfun}
%is well defined if \r{hwr} is satisfied.

To find a free field realization of the intertwining operators
we introduce  the generating functions
\bea
Z(z)&=&{:}\exp\left(-\tih \intt d\la\ \ee^{i\la z}{a'_\la\over
\sh(\tih\la)}\right){:}\ ,\nn%\label{Z}
\\
Z'(z)&=&{:}\exp\left(\tih \intt d\la\ \ee^{i\la z}{a_\la\over
\sh(\tih\la)}\right){:}\ .\nn%\label{Z'}
\eea

Now we are ready to prove the following
\smallskip

\noindent
{\bf Proposition 9.} {\it
The components of the twisted intertwining operators have the free field
realization}
\bea
\Psi^*_-(z)&=&Z(z)\ ,\label{Psi-}\\
\Psi^*_+(z) &=& \int_{C}
{du\over2\pi}\ \ee^{\pi\eta (z-u)}    \left[(q)^{1/2} E(u)
Z(z)+(q)^{-1/2}Z(z)E(u)\right]\ ,\label{Psi+}\\
\Psi_\nu(z)&=&\Psi^*_{-\nu}(z+i\tih)\ ,\quad \nu=\pm\ ,\label{Psidual}\\
\Phi_-(z)&=&Z'(z)\ ,\label{Phi-}\\
\Phi_+(z) &=& \int_{C'}
{du\over2\pi}\  \ee^{\pi\eta' (z-u)}  \left[(q')^{1/2}
Z'(z)F(u)+(q')^{-1/2}F(u)Z'(z)\right]\ , \label{Phi+}\\
\Phi^*_\ep(z)&=&\Phi_{-\ep}(z+i\tih)\ ,\quad \ep=\pm\ ,\label{Phidual}
\eea
{\it where the contours $C$ and $C'$ are the same as in} \r{2e} {\it and}
\r{2f}.
\smallskip

To prove  the proposition 6 we should
use the first terms in the comultiplication formulas
\r{com-e-fu}, \r{com-f-fu} and \r{comul-h} specified for the
operators
$x=e^\pm(u)$, $f^\pm(u)$ and $h^\pm(u)$ and  the action
of these generating functions on the elements
of the evaluation two-dimensional module \r{eval-0}--\r{eval-h1}.
The result is
the commutation relations of the components of the
intertwining operators with the generating functions of the algebra
$\Ael$. For example,  for the operators $\Psi^*(z)$ these
defining equations are:
\bea
h^\pm(u)\Psi_-^*(z)&=&
{\sh\,\pi\eta(u-z+i\tih/2\pm i\tih/4)
\over
\sh\,\pi\eta(u-z-i\tih/2\pm i\tih/4)}
\Psi_-^*(z)h^\pm(u)\ ,\label{51}\\
0&=& f^\pm(u)\Psi_-^*(z)+ \Psi_-^*(z)f^\pm(u)\ ,\label{52}\\
\sh\,i\pi\eta\tih \Psi^*_+(z)&=&
\sh\,\pi\eta(u-z-i\tih/2\pm i\tih/4)
e^\pm(u)\Psi_-^*(z)+\nn\\
&&\quad +\
\sh\,\pi\eta(u-z+i\tih/2\pm i\tih/4)
\Psi_-^*(z)e^\pm(u)\ .\label{53}
\eea
Because of \r{FDE} and \r{FDF} from \r{52} and \r{53}
follows that the operator $\Psi^*_-(z)$
anticommute  with the generating function $F(u)$ and has the commutation
relation with $E(u)$  as follows:
$$%\beq
\sh\,\pi\eta(u-z-i\tih/2)
E(u)\Psi_-^*(z)=-
\sh\,\pi\eta(u-z+i\tih/2)
\Psi_-^*(z)E(u)\ .
$$%\label{54}\eeq
It is easy now to verify   using formulas of the Appendix B
that the generating function $Z(z)$ satisfy these commutation relations
with $E(u)$, $F(u)$ and also \r{51}.

The representation of $\Psi^*_+(z)$ in integral form follows from \r{53}
and \r{2e}. The analysis of the normal ordering relations
of the generating functions $E(u)$, $Z(z)$, $F(u)$ and $Z'(z)$ 
shows that the contours $C$ and $C'$ in \r{Psi+} and \r{Phi+} 
coincide with those in \r{2e} and \r{2f}. 

Comparing the
 formulas \r{Psi-}--\r{Phidual} with the free field representation
of the Za\-mo\-lod\-chi\-kov-Faddeev 
ope\-ra\-tors from \cite{JKM} we conclude
that these operators coincide with twisting intertwining operators.
Therefore, they satisfy the Zamolodchikov-Faddeev algebra:
\begin{eqnarray}
\Psi^*_{\nu_1}(z_1) \Psi^*_{\nu_2}(z_2)
&=&S_{\nu_1\nu_2}^{\nu'_1\nu'_2} (z_1-z_2,\xi)
\Psi^*_{\nu'_2}(z_2) \Psi^*_{\nu'_1}(z_1)\ , \label{ZFII} \\
\Phi_{\ep_2}(z_2) \Phi_{\ep_1}(z_1)
&=&\tilde R_{\ep_1\ep_2}^{\ep'_1\ep'_2} (z_1-z_2,\xi+\tih)
\Phi_{\ep'_1}(z_1) \Phi_{\ep'_2}(z_2)\ , \label{ZFI} \\
\Psi^*_{\nu}(z_1) \Phi_{\ep}(z_2) &=&\nu\ep\,\tg\left({\pi\over4}
+{i\pi(z_1-z_2)\over2\tih}\right)    \Phi_{\ep}(z_2)
\Psi^*_{\nu}(z_1)\ , \label{ZFI-II}\\
\label{orthI}
\Phi_{\ep_1}(z)\Phi^{*}_{\ep_2}(z)
&=&g'\delta_{\ep_1\ep_2}\ \id\ ,\\
\label{orthII}
\Psi_{\ep_1}(z_1)\Psi^{*}_{\ep_2}(z_2)
&=&{g\delta_{\ep_1\ep_2}\, \id
\over z_1-z_2} + o(z_1-z_2)\ ,
\end{eqnarray}
where  the $S$-matrix in \r{ZFII}  is given by \r{trans},
the $R$-matrix $\tilde R(z,\xi+\tih)$ in \r{ZFI}
is related to the $R$-matrix given by \r{R-mat} as follows:
\beq
\tilde R(z,\xi+\tih)=(\sigma_z\ot1) R(z,1/(\xi+\tih)) (1\ot\sigma_z)
\label{trans2}
\eeq
and the constants $g$, $g'$ are equal
to\footnote{The double $\Gamma$-function $\Gamma_2(x\mid\omega_1;\omega_2)$
is defined in the Appendix B.}
\bea
g&=&{i\ee^{-3\gamma\eta/2\eta'}\eta^{2\eta/\eta'}\over
\eta^2\Gamma^2(\eta/\eta')}
{\Gamma_2(2\tih\mid2\tih;1/\eta) \Gamma_2(2\tih+1/\eta\mid2\tih;1/\eta)
\over
\Gamma_2(\tih\mid2\tih;1/\eta)\Gamma_2(3\tih+1/\eta\mid2\tih;1/\eta)}\ ,\nn\\
g'&=&{i\ee^{-3\gamma\eta'/2\eta}{\eta'}^{2\eta'/\eta}\over
\sqrt{2\pi\eta'}\Gamma^2(\eta'/\eta)}
{\Gamma_2(\tih\mid2\tih;1/\eta') \Gamma_2(1/\eta'-\tih\mid2\tih;1/\eta')
\over
\Gamma_2(2\tih\mid2\tih;1/\eta') \Gamma_2(1/\eta'\mid2\tih;1/\eta')}\  .
\eea

The proof of the commutation relations \r{ZFII}--\r{orthII}
can be found in \cite{JKM} and 
 is based on the formulas gathered in the
Appendix B.
In order
to prove  \r{orthI} and \r{orthII} one should use following
operator identities \cite{L}
\beq
E(u)=\ee^\gamma\left({:}Z\left(u+{i\tih\over2}\right)
Z\left(u-{i\tih\over2}\right){:}\right)^{-1},\quad
F(u)=\ee^\gamma\left({:}Z'\left(u+{i\tih\over2}\right)
Z'\left(u-{i\tih\over2}\right){:}\right)^{-1}.\label{Luk}
\eeq
The identities \r{Luk} being substituted
into \r{Psi+} and \r{Phi+} yield the integral relation between
com\-po\-nents of the intertwining (Zamolodchikov-Faddeev) operators.
This relation can be treated as a quantum version 
of the relation between lenearly independent solutions 
to the 
second order ordinary differential equation $(\partial^2+u(z))\psi(z)=0$.
%If the function $\psi_1(z)$ is a solution to this equation than the
%second solution can be found  as $\psi_2(z)=-\psi_1(z)\int^z\psi_1^{-2}(z)dz$.

\subsection{Zero mode discussion}

In this subsection  we follow  a well known idea presented in 
\cite{S1,BL}
in order to discuss unusual (twisted) definition
of the intertwining operators \r{inteq}. The reason for this
definition lies in the absence of the zero mode operator $(-1)^p$
in the bosonization
of the algebra $\Ael$ and twisted intertwining operators. We can come
back to the usual definition  introducing additional operators
$P$ and $I$ such that
$$
[P,a_\la]=[I,a_\la]=0,\quad PI=iIP\quad \mbox{and}\quad I^4=1.
$$
Define also the extended Fock space
$$
\overline{\H}=\H\ot\CC[\ZZ/4\ZZ]=\H_0\oplus\H_1
$$
and the subspaces $\H_0$ and $\H_1$  formed by the elements
$$
\H_0=\CC [v\ot1] \oplus \CC [v\ot I^2],\quad
\H_1=\CC [v\ot I] \oplus \CC [v\ot I^3],\quad v\in\H.
$$
The generating functions of the currents and the intertwining operators
are modified as follows:
$$
\tilde E(u)=E(u)\cdot I^2,\quad \tilde F(u)=F(u)\cdot I^2,\quad \tilde
H^\pm(u)=H^\pm(u)\ ,
$$
$$
\tilde \Psi^*_\pm(z)=\Psi^*_\pm(z)\cdot I^{\pm1}P,\quad
\tilde \Phi_\pm(z)=\Phi_\pm(z)\cdot I^{\pm1}P\ .
$$
It is clear now that the commutation relations of the modified intertwining
operators with elements of the algebra will be usual, for example:
$$
\tilde\Phi(z)x =\Delta(x) \tilde\Phi(z),\quad x\in\Ael\ ,
$$
since  the action of the operator $P$ on the elements of the algebra
coincide with the action of the involution \r{auto}
$$
P\,x\,P^{-1}=P^{-1}\,x\,P=\iota (x),\quad \forall\ x\in\Ael\ .
$$
The subspaces $\H_0$ and $\H_1$ become irreducible with respect
to the action of the algebra and the operators $\Phi(z)$, etc.
intertwine these subspaces. Since the known physical models
for which the algebra $\Ael$ serves as the algebra of dynamical
symmetries have single vacuum states this mathematical construction
of two irreducible Fock spaces is unnecessary  and this 
unnecessity  explains the absence
of the zero mode operators in the  bosonization of the massive models
of the quantum field theory.

\subsection{Miki's Formulas}

We would like to demonstrate now that bosonized expressions for the
intertwining operators are in accordance with the $L$-operator
description of the algebra $\Ael$. The way to do it is to use
Miki's formulas \cite{M}. Consider the $2\times2$ operator valued
matrices acting in the Fock space $\H$  ($\ep,\nu=\pm$)
\beq
L^+_{\ep\nu}(u) = \nu\ep\,\sqrt{{2\tih\ee^\gamma\over\pi}}
\Psi^*_\nu\sk{u-\fract{i\tih}{4}} \Phi_\ep\sk{u-\fract{3i\tih}{4}},\quad
L^-_{\ep\nu}(u) = \sqrt{{2\tih\ee^\gamma\over\pi}}
\Phi_\ep\sk{u-\fract{i\tih}{4}} \Psi^*_\nu\sk{u-\fract{3i\tih}{4}}\ .
\label{Miki}
\eeq
Now it is easy to show that so defined $L$-operators satisfy the
commutation relations \r{RLL-univ-+-} if the operators $\Phi(u)$ and
$\Psi^*(u)$ satisfy the commutation relations of the
Zamolodchikov-Faddeev algebra \r{ZFII}--\r{ZFI-II}.
Let us note here that although the intertwining operators commute by means
of the  $S(z,\xi)$ and $\tilde R(z,\xi+\tih)$
matrices which differ from $R(z,\eta)$ (see \r{trans} and  \r{trans2})
the $L$-operators defined by the Miki's  prescription \r{Miki}
commute according to \r{RLL-univ-+-} defined by $R$-matrix \r{R-mat}.

The Miki's formulas can be interpreted also as bosonization of
$L$-operators for the algebra $\Ael$ at level 1. Using these formulas
we can easily verify the relation 
\r{L-oper-rel} between $L^\pm(u)$ operators.
To do this it is
sufficient to use the Gauss decomposition of these operators 
and calculate
$$
L^+_{--}(u) =
{:}\,\Psi^*_-\sk{u-\fract{i\tih}{4}} \Phi_-\sk{u-\fract{3i\tih}{4}}{:}
={:}\exp \left(\tih\intt d\la\ \ee^{i\la u}{a_\la\ee^{-\la/2\eta''+\la\tih/2}
\sh(\la\tih/2)\over\sh(\la\tih)\,\sh(\la/2\eta)} \right){:},
$$
$$
L^-_{--}(u) =
{:}\,\Phi_-\sk{u-\fract{i\tih}{4}} \Psi^*_-\sk{u-\fract{3i\tih}{4}}{:}
={:}\exp \left(\tih\intt d\la\ \ee^{i\la u}{a_\la\ee^{\la/2\eta''+\la\tih/2}
\sh(\la\tih/2)\over\sh(\la\tih)\,\sh(\la/2\eta)} \right){:}.
$$
Note that the constant $\sqrt{2\tih\ee^{\gamma}/\pi}$ in \r{Miki}
is cancelled after normal ordering of  the operators $\Psi^*(u)$
and $\Phi(u)$ due to \r{ZZ'}.
The equality  $L^+_{--}(u-i/\eta'')=L^-_{--}(u)$ is obvious now and
the rest of the relations between the elements of the $L$-operators
can be found using the definition  of ``$+$''-components of the
intertwining operators through the generating functions $e^\pm(u)$,
$f^\pm(u)$ of the algebra $\Ael$ (see \r{53}) and the commutation
relation of these generating functions with $k^\pm(u)$.

\section{Discussion}

To conclude  this paper 
we would like to mention some  open problems which, to our opinion, 
deserve further investigation.

Let us note first that there are many possibilities to choose initial
parameters $\eta =1/\xi$ and $\tih$ in the definition of the algebra 
$\Ael$. For instance, 
it is clear from above analysis that the 
 properties of the algebra $\Ael$ or its representations
change drastically when the parameter $\xi=1/\eta=r\tih$,
where $r$ is some rational number. The commutation relations \r{ef}--\r{hh}
 show that in this case a smaller factoralgebra of $\Ael$ could be defined.
   It is naturally to assume that this factoralgebra serves the symmetries of
$\Phi_{[1,3]}$-perturbations of the minimal models of conformal
field theories. 

It is known from the theory of the Sine-Gordon model that if the parameter
$\xi<\tih$ than the spectrum of the model possesses scalar particles,
so called breathers. Our considerations were based on the assumption that 
 $\xi>\tih$ and it is an open question to investigate the algebra $\Ael$
 in the regime  $\xi<\tih$ and to apply its representation theory to the 
 Sine-Gordon model in the breather's regime.
We  think that 
the algebra ${{\cal A}_{\tih,\nu/\tih}(\widehat{\frak{sl}_2})}$,
$\nu=1,2,\ldots$ is related to the 
restricted Sine-gordon model in the reflectionless points 
and representation theory of this algebra 
can be used
for the group-theoretical interpretation of the results
obtained recently in \cite{BBS}. In particular, it is interesting to 
investigate the simplest case $\xi=\tih$ which should correspond
to the free fermion point of the Sine-Gordon model.

It is also natural from algebraic point of view to put the value of 
 deformation parameter $\tih$ to be pure imaginary instead of positive real.
 It could correspond to the Sinh-Gordon theory. Surely, one can also try 
 to apply the known technique 
\cite{Mat,AOS}
for studying $c>1$ integer level integrable
 representations of the algebra $\Ael$ and to find out possible physical
 applications.

As we have shown  in the 
Section 2, the definition of the algebra $\Ael$ cannot be
 done in purely algebraical terms. So far the corresponding analytical
 language should be developed for its representations. It was partially 
 done in Section 6 for the level one Fock space. Nevertheless, the rigorous
 mathematical description of the space of representation is far from 
 completeness. One needs more detailed topological description of the space,
 the precise definition of the trace, 
 making the calculations in \cite{L,JKM,Ko} to be rigorous,
 the investigation of the 
 irreducibility and so on. Moreover,
it would be nice
 to have an axiomatical description of the category of highest weight 
 representations.

The analysis of the defining commutation relations of the algebra 
$\Ael$ demonstrate that the algebra in question possesses the rich 
structure of automorphisms. For example, if we allow parameters
$\la$, $\mu$ in \r{he-hf}--\r{hh-hh} to be complex then these
commutation relations can be rewritten in the form of difference
commutation relations without integrals in the r.h.s. 

 The algebra $\Ael$ is not a Hopf algebra and even for $c=0$, when
it becomes a Hopf algebra, it does not have a structure of the quantum
double \cite{Dr},
 as well as its classical counterpart
 \cite{KLPST}. 
 The double 
structure can be reconstructed in the Yangian limit $\eta\to0$ when
the algebra ${\cal A}_{\tih,0}(\widehat{\frak{sl}_2})$ becomes
the central extended Yangian double. The representation theory 
of the latter algebra  have been 
investigated in \cite{KLP1,KLP2,IK,Ko} using two alternative possibilities
related to the Riemann problems on the circle and on the line.
It is interesting to understand what structure replaces the double structure 
in 
 $\Ael$; whether there exists an analog of the universal $R$-matrix.

It is also interesting to formulate the quantum Sugawara construction
corresponding to the algebra $\Ael$. The scaling quantum Virasoro
algebra can be obtained from the papers \cite{FF,FR,Lu,SKAO}  where
two-parameter deformation of the Virasoro algebra corresponding
to the algebra $\Apq$ has been investigated.

\section{Acknowledgments}

The research described in this publication was made possible in part by
grants 
RFBR-96-01-01106 (S.~Khoroshkin),
RFBR-96-02-18046 (D.~Lebedev),
RFBR-96-02-19085 (S.~Pakuliak),
INTAS-93-0166-Ext (S.~Khoroshkin, D.~Lebedev),
INTAS-93-2058-Ext  (S.~Pakuliak)
and by Award  No. RM2-150 of the U.S. Civilian Research \& Development
Foundation (CRDF) for the Independent States of the Former Soviet Union.

\setcounter{equation}{0}
\app{Consistency of Comultiplication Formulas}

To prove that the comultiplication formulas \r{comul-L-univ}
are in accordance with the commutation relations of the algebra
$\Ael$  we have to
check that the relations  (here $L(u)=L^+(u)$)
\begin{eqnarray}
R^+(u_1-u_2,\xi+\tih\Delta c)\Delta L_1(u_1) \Delta L_2(u_2)&=&
\Delta L(u_2) \Delta L(u_1) R^+(u_1-u_2,\xi)  ,
\nn%\label{RLLcom}
\\
\Delta(\qdet L(u))&=&\Delta(\qdet L(u)) \otimes \Delta(\qdet L(u))
\nn%\label{qdet-com}
\eea
follow from \r{RLL-univ} and the definition of the quantum determinant
\r{qdet}.
First of all we rewrite the equation \r{RLL-univ} in components
 and the comultiplication formulas
\r{com-L*cmp} using the short notations.
We have
$$%\beq
R_{im;kp}(u_1-u_2,\xi+\tih c) L_{mj}(u_1,\xi) L_{pl}(u_2,\xi)
= L_{kq}(u_2,\xi)   L_{ir}(u_1,\xi) R_{rj;ql}(u_1-u_2,\xi)
$$%\label{RLL-cmp}\eeq
and
$$%\beq
\Delta L_{mj}(u,\xi)= L_{fj}(u',\xi')L_{mf}(u'',\xi'')\ ,
$$%\label{com-L-sh}\eeq
where, as it was in the subsection 3.3, prime (double prime) denotes
that the corresponding $L$-operator or its component is in the first
(second) component of the tensor product and
$u'=u+i\tih c''/4$, $\xi'=\xi=1/\eta$
($u''=u-i\tih c'/4$, $\xi''=\xi+\tih c'$).
The summation over repeating indeces
is always supposed.
\bea
&&
R_{im;kp}(u_1-u_2,\xi+\tih c'+\tih c'')
\Delta\left(L_{mj}(u_1,\xi)\right) \Delta\left(L_{pl}(u_2,\xi)\right)\nn\\
&&=
R_{im;kp}(u_1-u_2,\xi+\tih c'+\tih c'')
L_{fj}(u'_1,\xi)
            L_{mf}(u''_1,\xi+\tih c')
L_{f'l}(u'_2,\xi)
            L_{pf'}(u''_2,\xi+\tih c')\nn\\
&&=
R_{im;kp}(u_1-u_2,\xi+\tih c'+\tih c'')
L_{fj}(u'_1,\xi)
L_{f'l}(u'_2,\xi)
            L_{mf}(u''_1,\xi+\tih c')
            L_{pf'}(u''_2,\xi+\tih c')\nn\\
&&=
R_{r'f;q'f'}(u_1-u_2,\xi+\tih c')
L_{fj}(u'_1,\xi)
L_{f'l}(u'_2,\xi)
            L_{kq'}(u''_2,\xi+\tih c')
            L_{ir'}(u''_1,\xi+\tih c')\nn\\
&&=
R_{rj;ql}(u_1-u_2,\xi)
L_{q'q}(u'_2,\xi)
L_{r'r}(u'_1,\xi)
            L_{kq'}(u''_2,\xi+\tih c')
            L_{ir'}(u''_1,\xi+\tih c')\nn\\
&&=
R_{rj;ql}(u_1-u_2,\xi)
L_{q'q}(u'_2,\xi)
            L_{kq'}(u''_2,\xi+\tih c')
L_{r'r}(u'_1,\xi)
            L_{ir'}(u''_1,\xi+\tih c')\nn\\
&&=
\Delta\left(L_{kq}(u_2,\xi)\right) \Delta\left(L_{ir}(u_2,\xi)\right)
R_{rj;ql}(u_1-u_2,\xi)\ .
\nn
\eea
The  primitivness of the  coproduct of the quantum determinant can be proved
easily using   \r{ke-sp}, \r{kf-sp} and the formula equivalent to \r{qdet}
$$
\qdet\, L(u)=D(u-i\tih)A(u)-C(u-i\tih)B(u)\ .
$$

Now we are in position to prove the formulas \r{com-e-fu}--\r{comul-h}.
The simplest comultiplication relation which follows from
\r{comul-L-univ} is
\bea
\Delta k(u,\xi)&=&
k(u',\xi')\otimes k(u'',\xi'')+
f(u',\xi') k(u',\xi')\otimes k(u'',\xi'') e(u'',\xi'')\nn\\
&=&
\left[1\ot 1+f(u',\xi') \otimes  e(u''-i\tih,\xi'')\right]
k(u',\xi')\otimes k(u'',\xi'') \label{inv1}\\
&=&
k(u',\xi')\otimes k(u'',\xi'')
\left[1\ot 1+f(u'-i\tih,\xi') \otimes  e(u'',\xi'')\right]\ .
\label{inv2}
\eea
The equivalent form of comultiplication of the operators
$k^\pm(u)$  follows from the operator identities
\bea
k(u,\xi)
e(u,\xi)-
e(u-i\tih,\xi)
k(u,\xi)  &=&0\ ,\label{ke-sp}\\
k(u,\xi)
f(u-i\tih,\xi)-
f(u,\xi)
k(u,\xi)  &=&0\ ,\label{kf-sp}
\eea
which are consequences of the commutation relations \r{RLL-univ}
at the critical point $u_1-u_2=i\tih$.

Formulas \r{com-e-fu} and \r{com-f-fu} easily follows from
\r{com-L*cmp} and \r{inv1}, \r{inv2}. The proof of \r{comul-h}
is more involved. It follows from the chain of
identities:\footnote{Here and below we do not write explicitly
dependence of the operators on the parameter $\xi$.}
\bea
\Delta h(u)&=&\sum_{p,p'=0}^{\infty}
(-1)^{p+p'}f^{p+p'}(u'-i\tih)h(u')\ot k^{-1}(u'') e^{p}(u''-i\tih)
k^{-1}(u''+i\tih) e^{p'}(u'')\nn\\
&=&\sum_{p,p'=0}^{\infty} (-1)^{p+p'}
f^{p+p'}(u'-i\tih)h(u')\ot  h(u'') \left(
[2]_\eta e(u''-i\tih)-e(u'')\right)^p e^{p'}(u'')\nn\\
&=&\sum_{p=0}^{\infty}    (-1)^p
f^{p}(u'-i\tih)h(u')\ot h(u'') \sum_{k=0}^p \left(
[2]_\eta e(u''-i\tih)-e(u'')\right)^k e^{p-k}(u'')\nn\\
 &=&\sum_{p=0}^{\infty}  (-1)^p
[p+1]_\eta f^{p}(u'-i\tih)h(u')\ot h(u'') e^{p}(u''-i\tih)\ .\nn
\eea
Here we used the commutation relation
which follows from \r{he} 
$$
k(u+i\tih) e(u-i\tih) k^{-1}(u+i\tih) = [2]_\eta e(u-i\tih) -e(u)
$$
and the combinatorial identity
$$
\sum_{k=0}^p  ([2]_\eta
e(u-i\tih)-e(u))^k e^{p-k}(u)=[p+1]_\eta e^p(u-i\tih)
$$
which follows by induction from \r{ee-univ}.

There is another way to verify the concordance of the formulas
\r{com-L*cmp}.
The comultiplication formulas for the currents
follow from comultiplication of $L$-operator entries
$L_{12}(u,\xi)$,
$L_{21}(u,\xi)$,
$L_{22}(u,\xi)$   given by
\r{com-L*cmp}.
The essential part of this calculation is the comultiplication of
the inverse operator $\left(k(u)\right)^{-1}$. Using the formulas
\r{inv1} and \r{inv2}
we can obtain
\bea
\Delta \left(k(u,\xi)\right)^{-1}&=&
\sum_{p=0} ^{\infty} (-1)^p
\left( f(u'-i\tih,\xi')\right)^p
\left(k(u',\xi')\right)^{-1} \otimes\nn\\
&\otimes&
\left(k(u'',\xi'')\right)^{-1}
\left(e(u''-i\tih,\xi'')  \right)^p \ .
   \label{inv4}
\eea
The comultiplication of the
entry
$L_{11}(u,\xi)$
also defines the comultiplication of the operator
$\left(k(u)\right)^{-1}$
so we should prove that these two comultiplication formulas lead
to the same result.
After some simple algebra we have
\bea
\Delta k^{-1}(u+i\tih)&=&
\left( k^{-1}(u'+i\tih)+f(u')k(u')e(u')\right)
\otimes\nn\\
&&\quad \otimes\ 
\left( k^{-1}(u''+i\tih)+f(u'')k(u'')e(u'')\right)+\nn\\
&&\quad +\ k(u')e(u')\otimes f(u'')k(u'') -
\Delta f(u) \Delta k(u) \Delta e(u)
\nn\\
&=&
 \sum_{p=0} ^{\infty} (-1)^p
       \left( f(u') \right)^p k^{-1}(u'+i\tih)\otimes
 k^{-1}(u''+i\tih) \left (e(u'')  \right)^p \ . \nn
\eea
The last line obviously coincide with \r{inv4} after shifting
$z\to z-i\tih$.

\setcounter{equation}{0}
\app{Normal Ordering Relations}

The relations below are based on the formulas which can be found in
\cite{Bar,JM}
\bea
\int_{\tilde C}{d\la\,\ln(-\la)\over2\pi i\la}\
{\ee^{-x\la}\over1-\ee^{-\la/\eta}}&=&\ln\Gamma(\eta x)+
\left(\eta x-\fract{1}{2}\right)(\gamma-\ln\eta)-\fract{1}{2}\ln2\pi\ ,
\nn%\label{integral}    
\\
\int_{\tilde C}{d\la\,\ln(-\la)\over2\pi i\la}\
{\ee^{-x\la}\over(1-\ee^{-\la\omega_1})
(1-\ee^{-\la\omega_2})}&=&\ln\Gamma_2(x\mid\omega_1,\omega_2)-{\gamma\over2}
B_{2,2}(x\mid\omega_1;\omega_2)\ ,
\nn%\label{integral2}
\eea
where $B_{2,2}(x\mid\omega_1;\omega_2)$ is the double Bernulli polynomial
of the second order
$$
B_{2,2}(x\mid\omega_1;\omega_2) ={1\over\omega_1;\omega_2}
\left[x^2-x(\omega_1+\omega_2)+{\omega^2_1+3\omega_1\omega_2+\omega^2_2
\over6}\right].
$$

Using these integral representations of the ordinary and double
$\Gamma$-functions and the definition of the product \r{normal}
we can  calculate:
\bea
%E(u)E(v)&=&i\eta(v-u)
%{\Gamma(1+\eta\tih-i\eta(u-v))
%\over \Gamma(-\eta\tih-i\eta(u-v))}
%{{:}\,E(u)E(v){:}\over
%\left({\ee^{-\gamma}\eta}\right)^{2\eta/\eta'}},\quad \Im(u-v)>\tih\ ,
%\nn%\label{EE}
%\\
%F(u)F(v)&=&i\eta'(v-u)
% {\Gamma(1-\eta'\tih-i\eta'(u-v))
%\over \Gamma(\eta'\tih-i\eta'(u-v))}
%{{:}\,F(u)F(v){:}\over
%\left({\ee^{-\gamma}\eta'}\right)^{2\eta'/\eta}},\quad \Im(u-v)<-\tih\ ,
%\nn%\label{FF}
%\\
Z(z)E(u)&=&
{\Ga{i\eta(u-z)-\eta\tih/2}
\over \Ga{1+i\eta(u-z)+\eta\tih/2}
}\ {{:}Z(z)E(u){:}\over\sk{\ee^{\gamma}/\eta}^{\eta/\eta'}}
,\quad \Im(u-z)<-{\tih\over2}\ ,
\nn%\label{ZE}
\\
E(u)Z(z)&=&
{\Ga{-i\eta(u-z)-\eta\tih/2}
\over \Ga{1-i\eta(u-z)+\eta\tih/2}
}\ {{:}Z(z)E(u){:}\over\sk{\ee^{\gamma}/\eta}^{\eta/\eta'}}
,\quad \Im(u-z)>{\tih\over2}\ ,
\nn%\label{EZ}
\\
Z(z)F(u)
&=&i\ee^\gamma(u-z){:}Z(z)F(u){:},\quad \Im(u-z)<0\ ,
\nn%\label{ZF}
\\
F(u)Z(z)
&=&-i\ee^\gamma(u-z){:}F(u)Z(z){:},\quad \Im(u-z)>0\ ,
\nn%\label{FZ}
\\
Z'(z)F(u)&=&
{\Ga{i\eta'(u-z)+\eta'\tih/2}
\over \Ga{1+i\eta'(u-z)-\eta'\tih/2}
}\ {{:}Z'(z)F(u){:}\over \sk{\ee^{\gamma}/\eta'}^{\eta'/\eta}}
,\quad \Im(u-z)<{\tih\over2}\ ,
\nn%\label{Z'F}
\\
F(u)Z'(z)&=&
{\Ga{-i\eta'(u-z)+\eta'\tih/2}
\over \Ga{1-i\eta'(u-z)-\eta'\tih/2}
}\ { {:}Z'(z)F(u){:}\over \sk{\ee^{\gamma}/\eta'}^{\eta'/\eta}}
,\quad \Im(u-z)>-{\tih\over2}\ ,
\nn%\label{FZ'}
\\
Z'(z)E(u)
&=&i\ee^\gamma(u-z){:}Z'(z)E(u){:},\quad \Im(u-z)<0\ ,
\nn%\label{Z'E}
\\
E(u)Z'(z)
&=&-i\ee^\gamma(u-z){:}E(u)Z'(z){:},\quad \Im(u-z)>0\ ,
\nn%\label{EZ'}
\\
Z(z_1)Z'(z_2)&=&{1\over\sqrt{2\tih\ee^\gamma}}
{\Gamma\left(\fract{1}{4}+\fract{i(z_2-z_1)}{2\tih}\right)
\over
\Gamma\left(\fract{3}{4}+\fract{i(z_2-z_1)}{2\tih}\right) }
{:}Z(z_1)Z'(z_2){:},\quad \Im(z_1-z_2)>-{\tih\over2}\ , \label{ZZ'}\\
Z(z_1)Z(z_2)&=&g(z_1-z_2){:}Z(z_1)Z(z_2){:},\quad \Im(z_1-z_2)>0\ ,
\nn%\label{ZZ}
\\
Z'(z_1)Z'(z_2)
&=&g'(z_1-z_2){:}Z'(z_1)Z'(z_2){:},\quad \Im(z_1-z_2)>-\tih\ ,
\nn%\label{Z'Z'}
\eea
where the functions $g(z)$ and $g'(z)$ are
\bea
g(z)&=&\exp \left(-\int_{\tilde C}{d\la\,\ln(-\la)\over2\pi i\la}
{(1-\ee^{-\la(\xi+\tih)})(1-\ee^{-\tih\la})\over
(1-\ee^{-\xi\la})(1-\ee^{-2\tih\la})}\ \ee^{i\la z}\right),\nn\\
&=&\ee^{\gamma\eta/2\eta'}
{\Gamma_2(\tih-iz\mid2\tih;1/\eta) \Gamma_2(\tih+1/\eta-iz\mid2\tih;1/\eta)
\over
\Gamma_2(-iz\mid2\tih;1/\eta) \Gamma_2(2\tih+1/\eta-iz\mid2\tih;1/\eta)}\ ,
\nn\\
g'(z)&=&\int_{\tilde C}{d\la\,\ln(-\la)\over2\pi i\la}
{(1-\ee^{-\xi\la})(\ee^{-2\tih\la}-\ee^{-\tih\la})\over
(1-\ee^{-(\xi+\tih)\la})(1-\ee^{-2\tih\la})}\ \ee^{i\la z},\nn\\
&=&\ee^{\gamma\eta'/2\eta}
{\Gamma_2(2\tih-iz\mid2\tih;1/\eta') \Gamma_2(1/\eta'-iz\mid2\tih;1/\eta')
\over
\Gamma_2(\tih-iz\mid2\tih;1/\eta')
\Gamma_2(\tih+1/\eta'-iz\mid2\tih;1/\eta')}\ .
\nn
\eea
As usually the normal ordering of all operators is calculated in the
regions specified above and then analytically continued to all
possible values of the spectral parameters.

\end{document}